\documentclass[epsf]{article}
\usepackage{amsfonts, euscript}
\usepackage{subfigure} 
\usepackage{epsfig}
\usepackage{amsmath, amsthm}
\newcommand{\be}{\begin{eqnarray}}     	\newcommand{\ee}{\end{eqnarray}}   
\newcommand{\rf}[1]{~(\ref{#1})}        
\newcommand{\lb}[1]{\label{#1}}

\newtheorem{dummy}{Dummy}

\theoremstyle{remark}

\theoremstyle{definition}
\newtheorem{corollary}{Corollary}
\newtheorem{lemma}{Lemma}
\newtheorem{theorem}{Theorem}
\newtheorem{proposition}{Proposition}

\newtheorem{definition}[dummy]{Definition}

\begin{document}
  \title{Option Pricing without Price Dynamics: A Probabilistic Approach}
\author{   Dimitris Bertsimas \thanks{Boeing Professor of Operations Research,
Sloan School of Management, Rm. E53-363, Massachusetts Institute
of Technology, Cambridge, Mass. 02139. dbertsim@mit.edu. Research
partially supported by   the
MIT-Singapore Alliance.} \and  Natasha Bushueva \thanks{
Work was done while being a graduate student at the Mathematics department of Massachusetts Institute of Technology  and is part of her PhD thesis (Apr 2003; http://hdl.handle.net/1721.1/29349). natacha@alum.mit.edu. 
} }
\date{August, 2004}
\date{}
\maketitle

\begin{abstract}
Employing probabilistic techniques we compute best possible 
 upper and lower bounds on the price of an option on one or two  assets with continuous piecewise linear payoff function  
based on prices of simple call options of possibly distinct maturities and the no-arbitrage condition, but without any assumption on the price dynamics of underlying assets.
We show that the problem reduces to solving linear optimization problems
that we explicitly characterize. 
We report numerical results that illustrate the effectiveness of the algorithms
 we develop.

\end{abstract}

\thispagestyle{empty}
\vfill

\section{Introduction}

One of the central questions in financial economics is to find the price of a derivative security given information on the underlying asset. Pricing of a derivative is usually realized by solving the Black-Scholes equation \cite{C1}, which is  based on the assumptions of a Geometric Brownian motion for the price of the underlying asset and  no-arbitrage in the market. One of the parameters in the equation, the volatility of returns on the underlying asset, is assumed to be defined and constant, while in reality it is varying in time and not known in advance. Practitioners typically use the implied volatility to value a derivative security by the Black-Scholes formula, that is they extract the volatility from prices of other options in the market, based on the assumption that the Black-Scholes equation is valid and other options are correctly priced. 

A natural question that arises is to determine the range of values for the price of an option based only on prices of other options and the no-arbitrage assumption, but with no assumption on the price dynamics of the underlying asset. 

   
In their seminal work,  Cox and Ross \cite{C2} and 
Harrison and Kreps \cite{C3} show that the condition of
no-arbitrage is equivalent to the existence a probability measure $Q$,
equivalent to the original measure $P$, with respect to which all
discounted securities processes are martingales. Rubinstein \cite{C4} and Longstaff \cite{C5} introduce the idea of deducing the martingale measure from observed European call prices by solving a quadratic optimization problem. 

 Lo \cite{C6} derives best possible closed form bounds on the price of
a European call option given the mean and variance of the underlying
stock price under risk neutral measure. Grundy \cite{C7}  extends Lo's
work for the case when the first and the $k$th moments of the stock
price are known.  Bertsimas and Popescu \cite{C13}  derive best possible bounds of the price of a European call option,  as well as on moments of the prices of the asset, given prices of other similar options, using a convex optimization approach. 
D'Aspremont and El Ghaoui \cite{Asp}  address the problem of computing upper and lower bounds on the price of a European basket call option, given prices on other similar baskets. They introduce a linear programming relaxation of the problem and show that this relaxation is best possible in some special cases.
All the above problems are solved for the case of a {\bf single maturity} and  single underlying stock. 

In this paper,
 we consider the problem of determining upper and lower bounds on the price of a European option given prices of other options of {\bf different maturities.} 
In particular, the contributions of the article are the following:

\begin{itemize}
\item[{\bf (a)}]
We establish  necessary and sufficient conditions that European call options of different maturities should jointly satisfy, so that there exists no arbitrage.
From this structural property and using geometric arguments, 
one can  derive best bounds on the price of a European call option with a given strike price and maturity. 
As we present more general methods in later sections
we do not present this method here (see Bertsimas and Bushueva \cite{bb1}). 


\item[{\bf (b)}] Given a set of European call options of different strike prices and different maturities we
determine best possible upper and lower bounds on the price of an option with a  continuous piecewise linear payoff function.
Based on direct probabilistic methods, we solve the problem as
 a linear optimization problem. 
%


\item[{\bf (c)}] Given European call options on two individual assets we determine best possible upper and lower bounds on the price of an option with a payoff being a continuous piecewise linear function of prices of the two assets at the time of option's maturity.
Again, options are allowed to have different maturities.
We present two algorithms based on linear optimization. 
The first is an exact formulation 
and  involves a linear optimization problem of
 potentially exponential size; the second is asymptotically exact
and involves a linear optimization problem of polynomial size. 
\item[{\bf (d)}] We determine upper and lower bounds on the price of a European basket call option, given prices on other similar baskets solving exactly the problem addressed in D'Aspremont and El Ghaoui \cite{Asp} approximately. 
\end{itemize}
Throughout the paper, we refer to these problems as Problems (a), (b), (c), (d). 
The paper is structured as follows: 
in Sections 2, 3, 4, 5 we present our solutions for Problems (a), (b), (c)  and (d)
respectively.

\section{Characterization theorem for the case of multiple maturities}

In this section, we  determine best possible bounds on the price of a European
 call option
 given  prices of options on the same underlying security, but
 with  potentially different maturities.\\

\noindent
The price of a European call option  with a strike price $k$ and maturity $t$ is given by 

\begin{equation}
\label{eq:key-formula}
{ C= e^{\displaystyle{-\int_0^t r(s)ds}}E_Q\left[\left(S_t-k\right)^+\right]=E_Q\left[\left(X_t-ke^{\displaystyle{-\int_0^t r(s)ds}}\right)^+\right],}
\end{equation}
where 
$r(s)$ is the instantaneous riskless rate of return, 
$\{S_t\}_{t\geq 0}$  is the price process of the underlying security
 defined on some probability space $\left(\Omega,{\cal B}, P\right)$,
 and $Q$ is a measure equivalent to $P$, such that $\{X_t\}_{t\geq 0}:
= \left\{S_t e^{-\int_0^t r(s)ds}\right\}_{t\geq
0}$
 is a martingale process under $Q$.
 Two measures are equivalent if they have the same null sets.

Thus, the problem we address in this section  can be reformulated as
 follows:
 given the  set
$$ 
{\cal S}: = \left\{(k_{ij},C_{ij}, t_i)\in \mathbb
R_+^3\;|\;t_1<t_2<\ldots<t_n, i\in \{1,2,\ldots, n\},
\;j\in\{ 1,2,\ldots, U(i)\}\right\}$$
 find the sufficient and necessary condition that values of the set should jointly satisfy so that there exists a nonnegative martingale $\{X_{t_1}, X_{t_2}, \ldots,X_{t^*},\ldots,X_{t_n}\}$, such that 
\be
\lb{eq:op1}
E\left[\left(X_{t_i}-k_{ij}
e^{-\int_0^{t_{i}} r(s)ds}\right)^+\right] = C_{ij}
\ee
for all points of the set.

\noindent
{\bf Remarks:}
\begin{itemize}
\item[{\bf (a)}] We can set without loss of generality  
$r(s)\equiv0$, 
since we can define  $\tilde k_{ij}=k_{ij} e^{-\int_0^{t_i} r(s)ds}$, 
and rewrite Eq. (\ref{eq:op1})
 as $E\left[\left(X_{t_i}-\tilde k_{ij}\right)^+\right] = C_{ij}$.
\item[{\bf (b)}]
Since we can assume $r(s)\equiv 0$, only the order of
$t_1,t_2,\ldots,t_n$ is important,
 but not the specific values that the sequence takes. 
\item[{\bf (c)}]
The price of an option with a strike price $0$ is given by
$E[X_t]=X_0$,
 which is the current stock price and, thus,  is always known.
\end{itemize}

In light of the above remarks,  a simpler formulation of the problem
is as follows.
\\[2pt]

\noindent \parbox[c]{13.5cm}{

\noindent {\bf  Problem (a)}
\hrule
\noindent \rule{0pt}{5mm}
\noindent 
Given an ordered set
\be 
\lb{setS}
{\cal S}: = \left\{(k_{tj},C_{tj})\in \mathbb R^{2}_+ \;|\; t\in \{1,2,\ldots, n\},\;j\in\{ 1,2,\ldots, U(t)\}\right\}
\ee
and $X_0\in\mathbb R_+$
 find the sufficient and necessary condition that values of the set should jointly satisfy so that there exists a nonnegative martingale $\{X_1, \ldots,X_{t^*},\ldots,X_n\}$ such that
\be
\lb{oprice1}
E\left[\left(X_t-k_{tj}\right)^+\right] = C_{tj} \;\;\mbox{for all}\;\;
(k_{tj},C_{tj})\in {\cal S},
\ee
and
\be
E[X_t] = X_0\;\;\;{\mbox {for all  }} t = 1,\ldots,n.\nonumber
\ee
\hrule
\noindent \rule{0pt}{5mm}
}
\noindent 
\\~

Let ${\cal C}$ designate the class of probability laws with support in $\mathbb R_+$ and finite expectations.
We associate with each 
element $\Pi$ of ${\cal C}$ the transform $\Psi_{\Pi}: \;\mathbb R_+ \rightarrow \mathbb R_+$, defined as
\be
\Psi_{\Pi}(t) = E_{\Pi}\left[(x-t)^+\right] = \int_0^{\infty}(x-t)^+d\Pi(x).
\ee
We next show that  $\Psi_{\Pi}(t)$ uniquely determines $\Pi$. 
\begin{theorem}
\begin{itemize}
\item[{\bf (a)}] The right derivative ${\Psi_{\Pi}}'(t+)$ of $\Psi_{\Pi}(t)$ exists for all $t\geq 0$ and 
\be
{\Psi_{\Pi}}'(t+)=-\Pi\left(\;(t,+\infty]\;\right).\nonumber
\ee
\item[{\bf (b)}]
 Let $\Pi_1, \;\Pi_2\in {\cal C}$ and for each $t\in \mathbb R_+$, ${\Psi_{\Pi_1}}(t)={\Psi_{\Pi_2}}(t)$. Then, $\Pi_1 = \Pi_2$.
\label{th_derivative}
\end{itemize}
\end{theorem}
\begin{proof}
\noindent{\bf (a)} Let $t\geq 0$ and $h>0$. We have
\be
\frac{\Psi_{\Pi}(t+h)-\Psi_{\Pi}(t)}{h}=\int_0^{\infty}{\frac{(x-(t+h))^+-(x-t)^+}{h}\;d\Pi(x)}=\int_0^{\infty}f_{t,h}(x)\;d\Pi(x),\nonumber
\ee
where the function $f_{t,h}(x)$ is defined as follows (see also Figure \ref{fth}): 
\be
f_{t,h}(x):=\frac{(x-(t+h))^+-(x-t)^+}{h}=\left\{
\begin{array}{lll}
0,\;\; & &x\leq t,\nonumber\\
{\displaystyle -\frac{x-t}{h}},\;\; & &t\leq x\leq t+h,\nonumber\\
-1,\;\; & &x\geq t+h.\nonumber
\end{array}
\right.
\ee
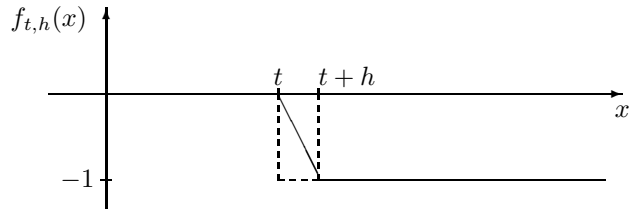
\begin{figure}[h]
    \begin{center}
      \setlength{\unitlength}{.3in}
      \begin{picture}(12,4)(-1.5,-2)
        \put(-1,0){\vector(1,0){10}}
        \put(0,-2){\vector(0,1){3.5}}
        \put(-0.5,-1.5){\makebox(0,0){$-1$}}
        \put(-0.1,-1.5){\line(1,0){.2}}
	\put(0,0){\line(1,0){3}}
	\put(3,0.1){\line(0,-1){0.2}}
        \put(3,0){\line(1,-2){0.72}}
	\put(3.7,-1.5){\line(1,0){5}}
	\put(3.7,0.1){\line(0,-1){0.2}}
        \put(3,0.3){\makebox(0,0){$t$}}
	\put(4.2,0.3){\makebox(0,0){$t+h$}}
        \put(-1,1.3){\makebox(0,0){$f_{t,h}(x)$}}
        \put(9,-0.3){\makebox(0,0){$x$}}
	\put(3,-1.5){\dashbox{0.1}(0.7,1.5){}}
      \end{picture}
    \end{center}
\caption{The function $f_{t,h}(x)$.}
\label{fth}
\end{figure}

\noindent Since the sequence of functions $f_{t,h}(x)$ is monotone (in $h$), and $f_{t,h}(x)$ converges as $h\rightarrow0^+$, then by 
 by the monotone convergence theorem \cite[p.100]{C8}, $\int_0^{+\infty}f_{t,h}(x)d\Pi(x)\rightarrow$ $-\Pi((t,+\infty))$ as $h\rightarrow0^+$. Thus $\Psi_{\Pi}'(t+)$ exists and equals  $-\Pi((t,+\infty)).$ 

\noindent{\bf (b)} Immediately follows from (a).
\end{proof}
The following proposition summarizes the properties of the $\Psi-$transform. 
\begin{proposition}
\label{pr_1}
Let $\Pi\in {\cal C}$. Then $\Psi_{\Pi}(t)$ satisfies:
\begin{itemize}
\item[{\bf (a)}] $\Psi_{\Pi}(t)$ is a nonincreasing function of $t\in \mathbb R_+$;
\item[{\bf (b)}] $\Psi_{\Pi}(t)$ is convex;
\item[{\bf (c)}] $\lim_{t\rightarrow+\infty}\Psi_{\Pi}(t)=0$. In particular, if $\Pi$ has a bounded support, then $\Psi_{\Pi}(t)=0$ for $t$ large;
\item[{\bf (d)}] $\Psi_{\Pi}(t)\geq \Psi_{\Pi}(0)-t$.
\end{itemize}
\end{proposition}
\begin{proof}

\noindent {\bf (a)} 
Notice, that ${\Psi_{\Pi}}(t)$ is a continuous function:
$$\left|\Psi_{\Pi}(t+h)-\Psi_{\Pi}(t)\right|=\left|\int_{t+h}^t(x-t)d\Pi(x)\right|\leq |h|.$$
Since $\Psi_{\Pi}(t)$ is continuous and ${\Psi_{\Pi}}'(t+)=-\Pi\left(\;(t,+\infty]\;\right)\leq 0$, then $\Psi_{\Pi}(t)$ is nonincreasing.
\\[1pt]
{\bf (b)} $f_x(t):=(x-t)^+$ is convex in $t$ for each $x$. Thus, $\Psi_{\Pi}(t)$ is convex as a convex combination of convex functions.
\\[5pt]
{\bf (c)}$$\Psi_{\Pi}(t)=\int_0^{+\infty}(x-t)^+d\Pi(x) = \int_t^{+\infty}(x-t)^+d\Pi(x)\leq\int_t^{+\infty}x d\Pi(x).$$
Since $\int_0^{+\infty}xd\Pi(x)<\infty$,
then $\int_t^{+\infty}x d\Pi(x)\rightarrow 0$ as $t\rightarrow \infty$. Thus,
$\Psi_{\Pi}(t)\rightarrow 0 $ as $t\rightarrow +\infty$.
\\[5pt]
{\bf (d)} $\Psi_{\Pi}(t) = \int_0^{\infty}(x-t)^+d\Pi(x)\geq \int_0^{\infty}(x-t)d\Pi(x) = \Psi_{\Pi}(0)-t$.
\end{proof}
We next show that  the properties in Proposition \ref{pr_1} 
are also sufficient for a function $\Psi$ to be a $\Psi-$transform of some distribution $\Pi$.  

\begin{theorem}
Let $g: \mathbb R_+\rightarrow \mathbb R_+$ be a function such that 
{\bf (a)} $g(x)$ is nonincreasing, 
{\bf (b)} $g(x)$ is convex,
{\bf (c)} $\lim_{x\rightarrow+\infty}g(x)=0$, and {\bf (d)}  $g(x)\geq g(0)-x$.
Then, there exists $\Pi\in {\cal C}$ such that $g(x) = \Psi_{\Pi}(x)$.
\label{theorem_g}
\end{theorem}
\begin{proof}
From the convexity of $g$ it follows that the right derivative of $g$, $g'(x+)$, is a right continuous nondecreasing function. 
Since $g(x)\geq g(0)-x$ and $g(x)$ is nonincreasing, then $-1\leq g'(x+)\leq0$. Also we have $g'(x+)\rightarrow 0$ as $x\rightarrow+\infty$. To justify the last statement, notice that $g'(x+)$ is a monotone bounded function, and consequently, it must converge to some finite value $a$ as $x\rightarrow\infty$. However, $a\neq 0$ would contradict the requirement in (c).

Consider the function $F(x): = 1+g'(x+)$ defined on $\mathbb R_+$. From the above statements, it follows that
({\em i}) $0\leq F(x)\leq 1$, ({\em ii}) $F(x)$ is a right-continuous
and nondecreasing, 
and ({\em iii}) $F(x)\rightarrow 1$ as $x\rightarrow \infty$.
Thus, $F(x)$ satisfies all the properties of a distribution function.
Denote by $\Pi$ the law with the distribution function $F(x)$. Then,
$$\Psi_{\Pi}(0) = E_{\Pi}[X] = \int_0^{+\infty}(1-F(x))dx=\int_0^{+\infty}\left(-g'(x+)\right)dx \stackrel{*}{=} g(0)-\lim_{x\rightarrow +\infty}g(x) = g(0).$$
The second to last equality (identified with ``*'') holds since $g(x)$
is a convex function and thus its right and left derivatives might
differ only on a countable set of points.

Since
${\Psi_{\Pi}}'(x+) = -(1-F(x)) = g'(x+)$,
$\Psi_{\Pi}(x)$ and $g(x)$ are both continuous and $\Psi_{\Pi}(0) = g(0)$, it follows that $\Psi_{\Pi}(x) = g(x)$ for any $x\in \mathbb R_+$. 
\end{proof}

Since there is a one to one correspondence between laws of random
variables and their $\Psi-$transforms, it follows that conditions that
a sequence of laws of random variables should satisfy can always be
reformulated in terms of conditions on the sequence of
$\Psi-$transforms of these random variables. Theorem \ref{th_big}
establishes  necessary and sufficient 
conditions that the sequence of $\Psi-$transforms should satisfy so that the corresponding random variables  form a martingale. However prior to proving that theorem, we will need the following auxiliary results.
\begin{definition}
\label{def_2}
Let $\Pi_1, \;\Pi_2\in {\cal C}$. Define the relation $\Pi_1 {\cal R}\Pi_2$ if and only if there exist random variables $X_1$ and $X_2$, defined on the same probability space, such that $\{X_1, X_2\}$ is a martingale and $X_1$ has a marginal law $\Pi_1$, while $X_2$ has a marginal law $\Pi_2$.
\end{definition}
\begin{proposition} The following properties hold:
\label{pr_2}
\begin{itemize}
\item[{\bf (a)}] If $\Pi_1{\cal R}\Pi_2$ and $\Pi_2{\cal R}\Pi_1$ then $\Pi_1 = \Pi_2$;
\item[{\bf (b)}] $\Pi_1{\cal R} \Pi_1$;
\item[{\bf (c)}] If $\Pi_1{\cal R}\Pi_2$ and $\Pi_2{\cal R}\Pi_3$, then $\Pi_1{\cal R}\Pi_3$.
\end{itemize}
\end{proposition}
\begin{proof}
 
\noindent {\bf (a)}  If $\{X_1,{\cal B}_1\}$, $\{X_2, {\cal B}_2\}$ is
a martingale with marginal laws $\Pi_1$ and $\Pi_2$ correspondingly,
and $f$ is a convex function, defined on $\mathbb R$, then
$f(X_1)=f(E\left[X_2|{\cal B}_1\right])\leq$ $E\left[f(X_2)|{\cal
B}_1\right]$ by the conditional Jensen's inequality, and thus $E[f(X_1)]\leq E[f(X_2)]$.  $(X-t)^+$ is a convex function of $X$, so $\Psi_{\Pi_1}(t)\leq \Psi_{\Pi_2}(t)$ for each $t\in \mathbb R_+$. From $\Pi_2{\cal R}\Pi_1$ by the same argument we obtain  $\Psi_{\Pi_2}(t)\leq \Psi_{\Pi_1}(t)$ for each $t\in \mathbb R_+$. Thus $\Psi_{\Pi_1}(t)=\Psi_{\Pi_2}(t)$ and by Theorem \ref{th_derivative}, $\Pi_1=\Pi_2$.
\\[2pt]
{\bf (b)} $\{X_1, X_1\}$ is always a martingale.
\\[2pt]
{\bf (c)} Let $\{X_1, X_2\}$ be a martingale with the law ${\cal L}(X_1, X_2)$ defined on $\mathbb R\times \mathbb R$ with support in  $\mathbb R_+\times\mathbb R_+$ and marginal laws $\Pi_1$ and $\Pi_2$ respectively. Let $\{\tilde{X_2},{X_3}\}$ be a martingale with the law ${\cal L}(\tilde{X_2}, {X_3})$ defined on $\mathbb R\times \mathbb R$ with support in $\mathbb R_+\times\mathbb R_+$ and marginal laws $\Pi_2$ and $\Pi_3$. Then on $\mathbb R$ there exist conditional distributions ${\cal L}(X_1|X_2)$ for ${\cal L}(X_1,X_2)$ and ${\cal L}(X_3|\tilde{X_2})$ for ${\cal L}(\tilde{X_2},X_3)$. Since the marginal laws of $X_2$ and $\tilde{X_2}$ are equal, then by Vorob'ev-Berkes-Philipp theorem \cite[p.7]{C12},  we can define a law ${\cal L}(X_1,X_2,X_3)$ on $\mathbb R\times\mathbb R\times \mathbb R$, such that $X_1$ and $X_3$ are conditionally independent given $X_2$, that is
$${\cal L}((X_1,X_3)|X_2)={\cal L}(X_1|X_2)\times{\cal L}(X_3|\tilde{X_2}).$$
 Denote by ${\cal B}_i$, $i=1,2,3$, the smallest $\sigma$-algebra for which all $X_k$, $k\leq i$, are measurable. Then $\{X_i,{\cal B}_i\}_{i=1}^3$ is a martingale with marginal laws $\Pi_1, \Pi_2$ and $\Pi_3$ correspondingly. Note that the martingale thus constructed is also a Markov process.
\end{proof}
%
%
%
%
%
%
%
\begin{theorem}
\label{th_big}
Let $\Pi_1,\Pi_2,\ldots,\Pi_n\in {\cal C}$ and $E_{\Pi_i}[X]=E_{\Pi_1}[X]$ for each $i=1,\ldots,n$. Then, there exists a martingale $\{X_1,X_2,\ldots,X_n\}$, such that $X_i$ has marginal law $\Pi_i$ for each $i$ if and only if the sequence of functions $\Psi_{\Pi_1}, \Psi_{\Pi_2},\ldots, \Psi_{\Pi_n}$ is nondecreasing.
\end{theorem}
\noindent {\bf Remark:} The 
theorem was obtained in Kertz and Rosler \cite{C11}. The key step in
their proof is a result of  Strassen \cite{C10}.
Our approach is completely different. 
\begin{proof}
The only if part is a direct consequence of the fact that $(X-t)^+$ is a convex function of $X$ and Jensen's inequality for convex functions of martingales \cite[p.277]{C8}. The proof of this statement is also in the proof of part {(a)} of Proposition \ref{pr_2}.

In the other direction, we have to show that it is possible to define a joint distribution of $X_1$, $X_2$,$\ldots$, $X_n$ with the given marginal laws (since they are uniquely determined by $\Psi_{\Pi_1}, \Psi_{\Pi_2},\ldots, \Psi_{\Pi_n}$), so that this sequence is a martingale. It is enough to show how to define a joint distribution of $X_1$ and $X_2$, since afterwards we can proceed recursively to form a martingale which is a Markov process as in the proof of Proposition \ref{pr_2}. 

We first  define a sequence $X_1, Y_1, Y_2,\ldots, Y_n,\ldots$, so
that this sequence is a martingale and $\Psi_{Y_n}\rightarrow
\Psi_{\Pi_2}$. We will then show that there exists $Y_{\infty}$, such that $\{Y_n\}_{1\leq n\leq\infty}$ is a martingale and $Y_{\infty}$ has a distribution $\Pi_2$.  Then $X_2$  can be defined as $X_2 = Y_{\infty}$.  

Let $Y_0=X_1$ and let us form a martingale $Y_0, Y_1, Y_2,\ldots, Y_n$
such that it is also a Markov process,
and thus,  we only need to define the joint distribution of $Y_i$ and $Y_{i+1}$ for each $i$. Now, the joint distribution of $Y_i$ and $Y_{i+1}$ is uniquely determined as soon as we know for each $y\in \mathbb R_+$ the conditional distribution of $Y_{i+1}$ given $Y_i = y$. Thus, we are going to define the conditional distribution for $Y_1$, then for $Y_2$ and so on.
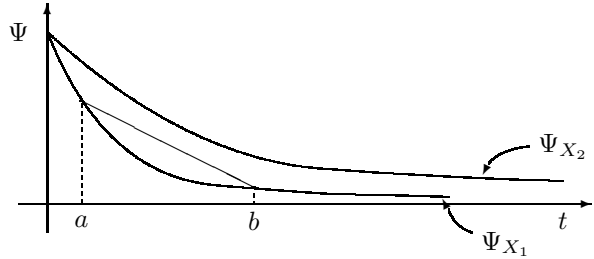
\begin{figure}[h]
    \begin{center}
      \setlength{\unitlength}{.3in}
      \begin{picture}(12,5)(-1.5,-0)
        \put(-0.5,0){\vector(1,0){10}}
        \put(0,-0.5){\vector(0,1){4}} 
	\qbezier(0,3)(2.5,0.7)(5,0.6) 
	\qbezier(5,0.6)(7,0.45)(9,0.4) 
    	\qbezier(0,3)(1,0.45)(3,0.32) 
	\qbezier(3,0.32)(4,0.23)(5,0.17) 
	\qbezier(5,0.17)(6,0.15)(7,0.12)
	\put(0.6,1.8){\line(2,-1){3}}
	\put(0.6,-0.3){\makebox(0,0){$a$}}
	\put(0.6,0){\dashbox{0.07}(0,1.8){}}
	\put(3.6,-0.3){\makebox(0,0){$b$}}
	\put(3.6,0){\dashbox{0.07}(0,0.3){}}
        \put(-0.5,3){\makebox(0,0){$\Psi$}}
        \put(9,-0.3){\makebox(0,0){$t$}}
	\put(9,1){\makebox(0,0){$\Psi_{X_2}$}}
	\qbezier(8.4,1)(7.9,0.95)(7.7,0.7)
	\put(7.7,0.7){\vector(-1,-2){0.1}}
	\put(8,-0.7){\makebox(0,0){$\Psi_{X_1}$}}
	\qbezier(7.4,-0.5)(7.1,-0.4)(7,-0.1)
	\put(7,-0.1){\vector(-1,2){0.1}}
      \end{picture}
    \end{center}
\caption{Construction of $\Psi_{Y_1}$.}
\label{Y_1}
\end{figure}
We know that $\Psi_{X_1} =\Psi_{\Pi_1}\leq\Psi_{\Pi_2}$. Now $Y_1$
will have its $\Psi$-transform, $\Psi_{Y_1}$, equal to $\Psi_{X_1}$
everywhere except on $(a,b)$ where it will be a straight line (Figure
\ref{Y_1}.) We choose $(a,b)$ so that the line which goes through the
points $\left(a,\Psi_{X_1}(a)\right), \left(b,\Psi_{X_1}(b)\right)$ is
below the function $\Psi_{\Pi_2}(t)$ (the specific choice of $(a,b)$ is discussed later.) What distribution will have its $\Psi$-transform equal to the transform of $\Pi_1$ everywhere except $(a,b)$ and being a straight line on $(a,b)$? From Theorem \ref{th_derivative} (a), we see that $\Psi_{Y_1}$ should have no weight on $(a,b)$ but the weight of $[a,b]$ should be the same as for the distribution $\Pi_1$. Thus, we take all the weight of $(a,b)$ and redistribute it between $\{a\}$ and $\{b\}$. We have to do it in such a way that $\{X_1, Y_1\}$ is a martingale.

The law ${\cal L}(X_1)$ is uniquely determined by $\Psi_{X_1}$ and the law ${\cal L}(Y_1)$ is uniquely determined by $\Psi_{Y_1}$. We define the law $P:={\cal L}(X_1, Y_1)$ with support in $\mathbb R_+\times \mathbb R_+$ in the following way. As the  marginal law of $X_1$ take $\Pi_{X_1}$. Then for each $X_1\in \mathbb R_+$ define the conditional distribution of $Y_1$ given $X_1$.  First, if $X_1\in \left\{\mathbb R_+\backslash(a,b)\right\}$, take $Y_1=X_1$. Now, if $X_1\in(a,b)$ we should have
$P(Y_1\in \{a\}\cup\{b\}|X_1\in(a,b))=1$ for $\Psi_{Y_1}(t)$ to be a straight line on $(a,b)$ and equal to $\Psi_{X_1}(t)$ everywhere else. Let us denote $\pi_{X_1}(X_1,Y_1):= X_1$. 

Thus, for any $c\in(a,b)$ we want the joint distribution of $X_1, Y_1$ to be such that the following holds
$$aP(Y_1=a|X_1=c)+bP(Y_1=b|X_1=c) = c,$$ 
$$P(Y_1=a|X_1=c)+P(Y_1=b|X_1=c)=1.$$                                           
Consequently,
$$P(Y_1=a|X_1=c) = \frac{b-c}{b-a} \;,\;\;\; P(Y_1=b|X_1=c)=\frac{c-a}{b-a}.$$ 
Thus, the conditional distribution of $Y_1$ given $X_1$ is defined and $\{X_1, Y_1\}$ is a martingale. In the same manner we can define the conditional distribution of $Y_{i+1}$ given $Y_{i}$, for $i=1,2,\ldots$. 

The next step is to show that $\Psi_{X_2}(t)$ can be approached by $\Psi_{Y_n}(t)$ with any given precision, that is that we can choose $Y_1,Y_2,\ldots,Y_n$ so that $\Psi_{Y_n}(t)\rightarrow\Psi_{X_2}(t)$. For each $i$, $\Psi_{Y_{i+1}}(t)$ is different from $\Psi_{Y_i}(t)$ only on a certain interval $(a_i,b_i)$, where $\Psi_{Y_{i+1}}(t)$ is a straight line. Thus, the choice of the sequence $\{Y_n\}$ is equivalent to the choice of the sequence $\{(a_n,b_n)\}$.
\begin{figure}[htbp]
\begin{center}
\input{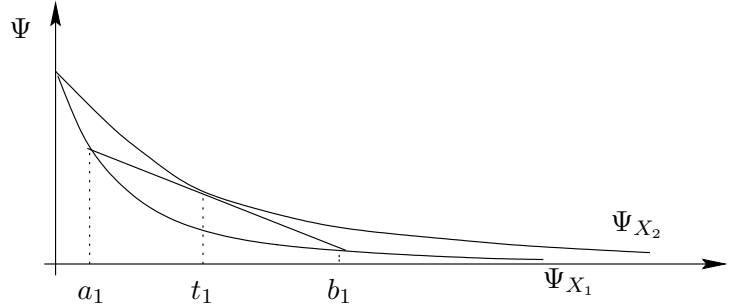}
\caption{The choice of $(a_1,b_1)$.}
\label{apr_12}
\end{center}
\end{figure}

\noindent Let 
$$U:=\left\{t\in \mathbb R_+: \Psi_{X_2}(t)>\Psi_{X_1}(t)\right\}.$$
Let $\{t_n\}$ be a countable dense set in $U$. 
At each $t_i$, $i=1,2,\ldots,$ draw a tangent line $y_i(t)$ to $\Psi_{X_2}(t)$. (See Figure \ref{apr_12}.)
Let us first assume that $\Psi_{X_2}(t_i)\neq 0$. Then ${\Psi_{X_2}}'(t_i+) < 0$ and the tangent line intersects $\Psi_{Y_{i-1}}(t)$ at two different points, $a_i\in \mathbb R_+$ and $b_i\in \mathbb R_+$,  since $\Psi_{Y_i}(t)\rightarrow0$ as $t\rightarrow0$ and $\Psi_{Y_i}(0)=1$, while $y_i(t)$ is below $\Psi_{X_2}(t)$ and thus $y_i(0)\leq1$. Thus, the choice of $\{t_n\}$ is equivalent to the choice of  $\{(a_n,b_n)\}$. We define $\Psi_{Y_{i+1}}=\max\left\{{y_i(t), \Psi_{Y_i}(t)}\right\}$. In case $\Psi_{X_2}(t_i)= 0$, ${\Psi_{X_2}}'_+(t_i)= 0$ and we take $\Psi_{Y_{i+1}}:=\Psi_{Y_i}$. 

Thus, the sequence $\left\{\Psi_{Y_n}\right\}$ is an increasing sequence of continuous functions converging to a continuous function $\Psi_{X_2}$ on a countable dense set. Consequently,  $\Psi_{Y_n}\rightarrow \Psi_{X_2}$ pointwise on $\mathbb R_+$.

Since $\Psi_{X_1}(t)\rightarrow 0$ and $\Psi_{X_2}\rightarrow0$ as $t\rightarrow +\infty$, then for any $\epsilon$ there exist $T\in \mathbb R_+$ such that for any $t>T$, $\;|\Psi_{X_2}(t)-\Psi_{X_1}(t)|<\epsilon$. Then for any $n$ and any $t>T$, $\;|\Psi_{X_2}(t)-\Psi_{Y_n}(t)|<\epsilon$. Taking into account that $\{\Psi_{Y_n}\}$ converges uniformly to $\Psi_{X_2}$ on $[0,T]$, we conclude that  $\Psi_{Y_n}$ converges uniformly to $\Psi_{X_2}$ on $\mathbb R_+$.  

Now, since $\Psi_{Y_n}\rightarrow \Psi_{X_2}$ uniformly on $\mathbb R_+$ and the right derivatives of $\Psi_{Y_n}(t)$, $n=1,2,\ldots,$ and $\Psi_{X_2}$ are monotone functions, then also ${\Psi_{Y_n}}'(t+)\rightarrow{\Psi_{X_2}}'(t+)$. Thus, the distribution functions of $Y_n$, $n=1,2,\ldots$, converge pointwise to the distribution function of $X_2$. 

We next show  that $\{Y_n\}$ is uniformly integrable, i.e.,  
 for any $\epsilon>0$ there exists $t_0>0$, such that for all $t>t_0$,
$$\sup_n E\left[Y_n1_{\{Y_n>t\}}\right]<\epsilon.$$
For each $t$ and $n$ we have
$\Psi_{Y_n}(t)\leq\Psi_{X_2}(t). $
Therefore,
$$\int_0^{+\infty}(Y_n-t)^+d\Pi_{Y_n}\leq \int_0^{+\infty}(X_2-t)^+d\Pi_{X_2},$$
$$\int_t^{+\infty}Y_n\;d\Pi_{Y_n}\leq \int_t^{+\infty}X_2\;d\Pi_{X_2}+ t\left|\Pi_{Y_n}\left([t,+\infty]\right)-\Pi_{X_2}\left([t,+\infty]\right)\right|.$$                                                                  
Since $\int_t^{+\infty}X_2\;d\Pi_{X_2}\rightarrow0$ as $t\rightarrow+\infty$, then we can take $t_0$ such that $\int_t^{+\infty}X_2\;d\Pi_{X_2}<\epsilon/3$ for all  $t>t_0$. 
Now since $\Pi_{Y_n}\left([t_0,+\infty]\right)\rightarrow\Pi_{X_2}\left([t_0,+\infty]\right)$, there exists $n_0$ such that for all $n>n_0$,  $\left|\Pi_{Y_n}\left([t_0,+\infty]\right)-\Pi_{X_2}\left([t_0,+\infty]\right)\right|<\frac{\epsilon}{3t_0}.$ Thus, for any $n>n_0$ and $t>t_0$ we have
$$E\left[Y_n1_{\{Y_n>t\}}\right]\leq E\left[Y_n1_{\{Y_n>t_0\}}\right]<\frac{\epsilon}3+t_0\frac{\epsilon}{3 t_0}=\frac{2\epsilon}3.$$
Therefore,
$$\sup_{n>n_0}E\left[Y_n1_{\{Y_n>t\}}\right]<\frac{2\epsilon}3.$$
For $Y_1, Y_2,\ldots, Y_{n_0}$ we can always choose $t_1$ big enough so that for all $t>t_1$, 
$E\left[Y_i1_{\{Y_i>t\}}\right]$ $<\frac{\epsilon}3$, $i=1,2,\ldots,n_0$, since all $Y_i$ have finite expectations. Consequently, if we denote $t^*:=\max\{t_0,t_1\}$, then 
$$\sup_{n}E\left[Y_n1_{\{Y_n>t\}}\right]<\epsilon,$$
for all $t>t^*$. 

Thus, $Y_1, Y_2, \ldots, Y_n$ is uniformly integrable, and,
 consequently \cite[p.283]{C8}, right closable,
which implies that there exists a random variable $Y_{\infty}$ such that $E[Y_{\infty}|{\cal B}_n]=$ $Y_n$ for all $n$, where ${\cal B}_n$ designates the smallest $\sigma$-algebra for which $Y_1,\ldots, Y_n$ are all measurable. Then by Doob's theorem \cite[p.285]{C8} $Y_n$ converges to $Y_{\infty}$ a.s., and consequently the distribution functions of $Y_1, Y_2,\ldots, Y_n$ also converge to the distribution function of $Y_{\infty}$. 
But then $Y_{\infty}$ and ${X_2}$  have the same distribution function, and thus the same $\Psi$-transform. So we take $X_2=Y_{\infty}$ and $X_1, X_2$ is a martingale.

\end{proof}

Now we are ready to formulate the conditions that 
the set ${\cal S}$ defined in (\ref{setS})
should satisfy, so that there exists a martingale
$\{X_1,X_2,\ldots,X_n\}$, such that Eq. (\ref{oprice1}) holds. 
\begin{proposition}
\label{pr_3}
The set ${\cal S}$, as defined in\rf{setS}, satisfies the condition of martingale existence if and only if there exists a nondecreasing sequence of convex functions, $\{g_t(x)\}_{t\in\{1,\ldots,n\}}$: $\mathbb R_+\rightarrow \mathbb R_+$, such that for each $t\in\{ 1,2,\ldots,n\}$
\\[1pt]
{\bf (a)} $g_t(0) = g_1(0)$;
\\[1pt]
{\bf (b)} $g_t(x)\geq g_t(0)-x$;
\\[1pt]
{\bf (c)} $g_t(x)\rightarrow 0$ as $x\rightarrow \infty$;
\\[1pt]
{\bf (d)} $g_t(k_{tj}) = C_{tj}$ for each $j\in \{1,2,\ldots,U(t)\}$.
\end{proposition}
\begin{proof}
$"\Longrightarrow"$ By Theorem \ref{theorem_g} for each $g_t$ there exists $\Pi_t\in {\cal C}$  such that $\Psi_{\Pi_t}=g_t$. Then, since the sequence of functions $\Psi_{\Pi_1}, \Psi_{\Pi_1}, \ldots, \Psi_{\Pi_n}$ is nondecreasing and $\Psi_{\Pi_t}(0) =$ $\Psi_{\Pi_1}(0)$ for each $t$, from Theorem \ref{th_big} it follows that the martingale satisfying all the points of the set ${\cal S}$  exists.   
\\[2pt]
$"\Longleftarrow"$ Suppose that there exists a martingale $\{X_1, X_2,\ldots,X_n\}$, such that\rf{oprice1} holds.  Then we can take $g_t(x): = \Psi_{X_t}(x)$ for $t = 1,2,\ldots,n$. By Theorem \ref{th_big} the sequence of functions $\{g_t\}_{t=1,\ldots,n}$ is nondecreasing and by Proposition \ref{pr_1} conditions (a)-(c) hold. By definition,  $g_t(k_{tj}) = \Psi_{X_t}(k_{tj})=C_{tj}$ for all $(k_{tj},C_{tj})\in{\cal S}$. 
\end{proof}
%
%
%
Designate by
\be
{\cal S}_t: = \left\{(k_{lj},C_{lj})\;|\; (k_{lj},C_{lj})\in{\cal S},\;l=t \right\}
\ee
the set of options of maturity $t$ and by
\be
{\cal S}_{\geq t}: = \cup_{l\geq t} {\cal S}_l
\ee
the set of options of maturities $t$ and higher. Let also
\be
{\cal S}_{\geq t}^{\infty}: = {\cal S}_{\geq t}\cup\{(0,+\infty)\}\cup\{(+\infty,0)\}.
\ee
Finally, let 
$\partial{\cal S}_{\geq t}:\mathbb R_+\longrightarrow \mathbb R_+$ represent the lower boundary of the convex hull of ${\cal S}_{\geq t}^{\infty}\cup\{(0,X_0)\}$, that is
\be
\partial{\cal S}_{\geq t}(x) = \min \left\{y|(x,y)\in\mbox{convex hull of } \{{\cal S}_{\geq t}^{\infty}\cup(0,X_0)\}\right\}.
\ee

\begin{theorem}
\label{hulls}
The set ${\cal S}$ satisfies the no-arbitrage condition if and only if 
\begin{itemize}
\item[{\bf (a)}] $C_{tj}\geq X_0-k_{tj}$ for each $(k_{tj},C_{tj})\in {\cal S}$;
\item[{\bf (b)}]  
for each $t$, none of the points of ${\cal S}_t$ is in the interior of the convex hull  of ${\cal S}_{\geq t}^{\infty}$.    
\end{itemize}
\end{theorem}

\begin{proof}
$"\Longrightarrow"$ Define functions $\tilde g_t(k): \mathbb R_+
\longrightarrow \mathbb R_+$ for each $t=1,\ldots,n$, such that $\tilde g_t(k) =  \partial {\cal S}_{\geq t}(k)\;\forall k\geq 0$. Then the sequence of functions $\{\tilde g_t\}$ satisfies all the properties of Proposition \ref{pr_3} except of property {\bf (c)}. Each function $\tilde g_t$ is piecewise linear with slopes of all, but the last, line segments being negative. The last line segment has a slope zero. Let $x_t: = \min\left\{x\;|\; x\in\mathbb R_+,\;\tilde g_t'(x)=0\right\}$. Denote by $\bold s$ the set of all slopes of functions $\tilde g_t$, $t=1,2,\ldots,n$, and  let $s^*:=\max \left\{s|s\in{\bold s},s<0\right\}$.  Define the function
\be
g_t(x):=\left\{
\begin{array}{rrr}
\tilde g_t(x), & &x\leq x_t,\nonumber\\
\tilde g_t(x_t)+s^*(x-x_t), & &\left\{x\geq x_t\;|\; g_t(x_t)+s^*(x-x_t)\geq 0\right\}, \nonumber\\
0, & & \mbox{otherwise}.\nonumber
\end{array}
\right.
\ee
Then the sequence of function $g_t(x)$ satisfies all the properties of Proposition \ref{pr_3} and, thus,  the martingale exists. 
\\[2pt]
$"\Longleftarrow"$ If $X_1,X_2,\ldots,X_n$ is a martingale, satisfying all the points of ${\cal S}$, then for any $(k,C)$ which belongs to the convex hull of ${\cal S}_{\geq 1}$, it must be that   $C\geq \Psi_{X_1}(k)\geq X_0-k$ by Proposition \ref{pr_1} and Theorem \ref{th_big}.
Assume that there exists a point $\left(k_{t^*j^*},C_{t^*j^*}\right)$, $t^*\in\{1,\ldots,n\}$, $j^*\in 1,\ldots,U\left(t^*\right)$, which is in the interior of the convex hull of ${\cal S}_{\geq t^*}^{\infty}$ (Figure \ref{duratskii}.) Take points $\left(k^l, C^l\right)\in{\cal S}_{\geq t}^{\infty}$, $\left(k^r, C^r\right)\in{\cal S}_{\geq t}^{\infty}$ such that $\partial {\cal S}_{\geq t^*}(k^l) = C^l$ and $\partial {\cal S}_{\geq t^*}(k^r) = C^r$  and $k^l$ is the closest to $k_{t^*j^*}$ from the left, while $k^r$ is the closest to $k_{t^*j^*}$ from the right.  Then $\left(k_{t^*j^*},C_{t^*j^*}\right)$ lies strictly above the interval connecting $\left(k^l,C^l\right)$ and $\left(k^r,C^r\right)$. 
Since $\Psi_{X_{t^*}}\leq\Psi_{X_{t^*+1}}\leq\ldots\leq \Psi_{X_n}$ by Theorem \ref{th_big}, then $\Psi_{X_{t^*}}(k^l)\leq C^l$ and $\Psi_{X_{t^*}}\left(k^r\right)\leq C^r$, thus $\left(k_{t^*j^*},C_{t^*j^*}\right)$ is also strictly above the interval connecting $\left(k^l,\Psi_{X_{t^*}}\left(k^l\right)\right)$ and $\left(k^r,\Psi_{X_{t^*}}\left(k^r\right)\right)$, which contradicts to the condition that $\Psi_{X_t^*}$ is a convex function. 
\begin{figure}[htbp]
\begin{center}
\input{dec18.pstex_t}
\caption{The point $\left(k_{t^*j^*},C_{t^*j^*}\right)$,
$t^*\in\{1,\ldots,n\}$, $j^*\in 1,\ldots,U\left(t^*\right)$ is 
 in the interior of the convex hull of ${\cal S}_{\geq t^*}^{\infty}$.}
\label{duratskii}
\end{center}
\end{figure}
\end{proof}

From the last theorem, we can find bounds on the price of a call option, given the set of priced options based on a geometric approach, see 
Bertsimas and Bushueva \cite{bb1}. 
However, we will suggest another approach, which leads to the linear optimization problem, and that is extendable to multiple dimensions.

\subsection*{Restriction to Markov martingales}

The idea of our approach is to reduce the class of martingale distributions over which  we maximize (or minimize) as much as possible, while making sure that the optimal solution does not change. In this section, we show that it is enough to consider martingales with the Markov property only. 

\begin{proposition}
\label{prop_6_2}
Let ${X_1}, {X_2},\ldots,{X_n}$ be a martingale with values in $\mathbb R^d$ and let $\mu_i$ be the marginal distribution of $X_i$. Then, there exists a martingale with the Markov property ${{\tilde X_1}}, {{\tilde X_2}}, \ldots,{{\tilde X_n}}$, such that the marginal distribution of ${{\tilde X_i}}$ is $\mu_i$ for all $i=1,2,\ldots,n$.
\end{proposition}
%
%
%

\begin{proof}
For specifying the distribution of a discrete time Markov process, it is enough to specify the transition probabilities as well as the initial distribution. 

Let $\nu ={\cal L}(X_1)$, $\mu_{1,2}={\cal L}(X_2|X_1), \ldots,\; \mu_{n-1,n}={\cal L}(X_n|X_{n-1})$. Then by the existence theorem for Markov processes due to Kolmogorov \cite[p. 120]{C9} there exists a Markov process ${\tilde X}_1,\ldots, \tilde X_n$ such that ${\cal L}(X_1)=\nu,\; {\cal L}(\tilde X_2|\tilde X_1)=\mu_{1,2},\;\ldots,\; {\cal L}(\tilde X_{n-1}|\tilde X_n)=\mu_{n-1,n}$.

It remains to show that the so-defined process $\tilde X_1, \tilde X_2,\ldots, \tilde X_n$ is also a martingale.
Let ${\sigma}_k$ be the smallest $\sigma$-algebra for which $\tilde X_1,\tilde X_2,\ldots, \tilde X_k$ are all measurable. We need to prove that $E[\tilde X_{k+1}|\sigma_k]=\tilde X_k$.
From the Markov property of $\tilde X_1,\ldots, \tilde X_k$, it follows that
$E[\tilde X_{k+1}|\sigma_k]=E[\tilde X_{k+1}|\tilde X_k]$. By definition $E[\tilde X_{k+1}|\tilde X_k]=E[X_{k+1}|X_k]$. 

 On the other hand, since $X_1,X_2,\ldots,X_n$ is a martingale, we have
\\[1pt]
$E\left[E[X_{k+1}|X_k, X_{k-1},\ldots,X_1]|X_k\right]=E[X_k|X_k]=X_k$. Also 
\\[1pt]
$E\left[E[X_{k+1}|X_k, X_{k-1},\ldots,X_1]|X_k\right]=E[X_{k+1}|X_k]$. Thus $E[X_{k+1}|X_k]=X_k$. Consequently, $E[\tilde X_{k+1}|\tilde X_k]=\tilde X_k$ and $E[\tilde X_{k+1}|\sigma_k]=\tilde X_k$. 
\end{proof}

\begin{definition}
\label{def_10}
A stochastic process which is a martingale and a Markov process will be called a Markov martingale.
\end{definition}
\begin{corollary}
\label{cor_4}
If we restrict our consideration to Markov martingales only, the optimal solution to our problems does not change.
\end{corollary} 

\section{The One Dimensional Case}

In this section, we address Problem (b), i.e., we 
find best possible upper and lower bounds on the price of a European style option with a continuous piecewise linear payoff function, given the set of European call options of different maturities.

The price of a European call option  with a strike price $k$ and maturity $t$ is given by Eq. (\ref{eq:key-formula}). The problem we address in this section 
can be formulated as follows: given a continuous piecewise linear function $g:\mathbb R_+\rightarrow\mathbb R_+$ and a set
\be 
{\cal S}: = \left\{(k_{ij},C_{ij}, t_i)\;|\; \left\{k_{ij}, C_{ij}, t_i\right\}\in \mathbb R_+,\;t_1<t_2\ldots<t_n,\;\right.\nonumber\\
\left. i\in \{1,2,\ldots, n\},\;j\in\{ 1,2,\ldots, U(i)\}\right\}\nonumber
\ee
find the domain of values of
$$e^{\displaystyle{-\int_0^{t^*} r(s)ds}}E\left[g\left(e^{{\int_0^{t^*} r(s)ds}}X_{t^*}\right)\right]$$
 under the condition, that $\{X_{t_1}, X_{t_2}, \ldots,X_{t^*},\ldots,X_{t_n}\}$
is a nonnegative martingale such that
\be
\lb{oprice10}
E\left[\left(X_{t_i}-k_{ij}e^{\displaystyle{-\int_0^{t_i} r(s)ds}}\right)^+\right] = C_{ij}.\nonumber
\ee
It is assumed that such martingale exists, that is the set of given options satisfies the no-arbitrage condition.\\
{\bf Remarks:}
\begin{itemize}
\item[{\bf (a)}]
As before, we can assume $r(s)\equiv0$, since we can take $\tilde k_{ij}=k_{ij}e^{-\int_0^{t_i} r(s)ds}$ instead of $k_{ij}$, and rewrite all the equalities as $E\left[\left(X_{t_i}-\tilde k_{ij}\right)^+\right] = C_{ij}$. 
Also, we can take the function $$\tilde g(X_{t^*}):=e^{-\int_0^{t^*} r(s)ds}g\left(e^{\int_0^{t^*} r(s)ds}X_{t^*}\right)$$ instead of $g(X_{t^*})$. Notice that if $g(X_{t^*})$ is a continuous piecewise linear function, so is $\tilde g(X_{t*})$. 
\item[{\bf (b)}]
Thus, again, only the order of $t_1,t_2,\ldots,t_n$ is important, but not the specific values that the sequence takes. Thus we can assume $\{t_1,t_2,\ldots,t_n\}\equiv \{1,2,\ldots,n\}$ for convenience. 
\item[{\bf (c)}] Recall that we must have $E[X_t]=X_0$, where $X_0$ is the current price of the stock.
\end{itemize}
In light of the above remarks, Problem (b) can be reformulated as follows:
\\[2pt]

\noindent \parbox[c]{13.5cm}{
\noindent {\bf  Problem (b)}

\noindent
\hrule
\noindent \rule{0pt}{5mm}
\noindent Given a continuous piecewise linear function $g:\mathbb R_+\rightarrow\mathbb R_+$ and an ordered set
\be 
\lb{setS2}
{\cal S}: = \left\{(k_{tj},C_{tj})\;|\; \left(k_{tj}, C_{tj}\right)\in \mathbb R^{2}_+,\; t\in \{1,2,\ldots, n\},\;j\in\{ 1,2,\ldots, U(t)\}\right\}
\ee
find the domain of values of
\be
\lb{terpi}
E\left[g(X_{t^*})\right]
\ee
for some  $t^*\in \{1,2,\ldots,n\}$, under the condition, that $\{X_1, \ldots,X_{t^*},\ldots,X_n\}$
is a nonnegative martingale such that
\be
\lb{oprice10}
E\left[\left(X_t-k_{tj}\right)^+\right] = C_{tj} \;\;\mbox{for all}\;\;
(k_{tj},C_{tj})\in {\cal S},
\ee
and
\be
E[X_t] = X_0\;\;\;{\mbox {for all  }} t = 1,\ldots,n.\nonumber
\ee
\hrule
\noindent \rule{0pt}{5mm}
}

\subsection{Treatment of future conditions}

Let
\be
{\cal F} =\left\{\left(k_{tj},C_{tj}\right)|t^*<t\leq n,\;\left(k_{tj},C_{tj}\right)\;
{\mbox{is a vertex of the convex hull of}}\;{\cal S}_{\geq t^*}^{\infty}\right\}\nonumber
\ee
 represent future points which are on the border of the convex hull of ${\cal S}_{\geq t^*}$.
And let $\left(k_{fut,j},C_{fut,j}\right)$, $j=1,2,\ldots,U(_{fut})$, be
the  elements of ${\cal F}$. Then we have the following lemma:

\begin{lemma}
\label{l_1}
If we reduce conditions\rf{oprice10} in Problem (b) to the set of conditions
\be
E\left[\left(X_t-k_{tj}\right)^+\right] = C_{tj},\quad t\leq t^*,\; j=1,2\ldots, U(t),
\lb{new_cond10}
\ee
\be
E\left[\left(X_{t^*}-k_{fut,j}\right)^+\right] \leq C_{fut,j},\quad j=1,2,\ldots, U(_{fut}),
\lb{new_cond20}
\ee
%
%
the solution to Problem (b) will not change.
\end{lemma}

\begin{proof} 
Suppose that we want to maximize $E\left[g(X_{t^*})\right]$ and $g^*_1$ is the optimal solution under conditions\rf{oprice10} and $g^*_2$ is the optimal solution under conditions \rf{new_cond10} and\rf{new_cond20}. Then $g^*_1\leq g^*_2$. Let us show that $g^*_1\geq g^*_2$.

Let $\Pi^*$ be the distribution that maximizes $\int_0^{+\infty}g(x)d\Pi(x)$ over distributions $\Pi$ of $X_{t^*}$ under  conditions\rf{new_cond10}\&\rf{new_cond20}.  Then there exists a martingale $X_1,\ldots,X_{t^*}$ which satisfies the part of conditions\rf{oprice10}  for $t\leq t^*$ and with the distribution $\Pi^*$ of  $X_{t^*}$, such that $\Psi_{\Pi^*}(k_{fut,j})\leq C(k_{fut,j})$ for all $j=1,\ldots,U(_{fut})$.  Since $\Psi_{\Pi^*}$ has the properties stated in  Proposition \ref{pr_1}, none of the points $(k,\Psi_{\Pi^*}(k))$, $k\geq 0$, is in the interior of the convex hull of ${\cal S}_{\geq t^*}^{\infty}$ and $\Psi(k)\leq X_0-k$.  Then, since also the set $S_{\geq t^*}$ by itself satisfies the no-arbitrage condition, by Theorem \ref{hulls} there exists a martingale  $\tilde X_t^*,\ldots,\tilde X_{t_n}$ which satisfies the set $\left\{S_{\geq t^*}\right\}\cup\left\{\left(k,\Psi_{\Pi^*}(k)\right)\right\}_{k\geq 0}$. Notice that the distribution of $\tilde X_{t^*}$ is!
  $\Psi_{\Pi^*}$. 
As in the proof of Proposition  \ref{prop_6_2}, it can be shown that the martingale  $X_1,X_2,\ldots,X_{t^*}$ can be extended to the martingale $X_1,\ldots,X_{t^*},\ldots,X_n$ by specifying transition probabilities  ${\cal L}(X_{t^*+1}|X_{t^*}) = {\cal L}(\tilde X_{t^*+1}|\tilde X_{t^*}), \ldots, {\cal L}(X_n|X_{n-1}) = {\cal L}(\tilde X_n|\tilde X_{n-1})$. 
This extended martingale satisfies all the conditions\rf{oprice10}. Thus $g^*_1\geq g^*_2$.
\end{proof}

\subsection{Restriction to discrete distributions}

From now on we make an additional assumption that 
\be
P(X_i\leq L, i=1,\ldots, t^*)=1.\nonumber
\lb{bound_on_stock}
\ee
for some $L>0$.
%
Let 
$${\cal K}_S=\{0\}\cup\{L\}\cup\left\{k_{tj}\;|\; \left(k_{tj},C_{tj}\right)\in {\cal S}, t\leq t^*\right\}\cup\left\{k_{tj}\;|\; (k_{tj},C_{tj})\in {\cal F}\right\}$$
and 
$${\cal K}_g=\left\{k\in[0,L]\;|\;g'(k-)\neq g'(k+)\right\}$$  
be the set of points, where $g$ changes its derivative. Define:
\be
{\cal K} = {\cal K}_S\cup{\cal K}_g.\nonumber
\ee

\begin{theorem}
\label{th_1}
Let $X_1,\ldots,X_{t^*}$ be a nonnegative Markov martingale with support in $[0,L]$.
Then there exists a martingale ${\tilde X}_1, \ldots, {\tilde X}_{t^*}$ with support in ${\cal K}$, such that
$$
E\left[\left(X_t-k_{tj}\right)^+\right] =E\left[\left(\tilde X_t-k_{tj}\right)^+\right],\quad t=1,2,\ldots, t^*;\quad j=1,2,\ldots,U(t) 
$$
$$
E\left[\left(X_{t^*}-k_{fut,j}\right)^+\right]=E\left[\left(\tilde X_{t^*}-k_{fut,j}\right)^+\right], \quad j=1,2,\ldots,U(_{fut})$$
and
$$E\left[g(X_{t^*})\right]=E\left[g(\tilde X_{t^*})\right].$$

\end{theorem}

\noindent We will use the following two lemmas to prove Theorem \ref{th_1}.

\begin{lemma}
\label{lemma_9}
Let $X_1, X_2$ be a martingale with values in $\mathbb R^d$. Then there exists a continuous time martingale $\left\{\tilde X_t\right\}_{t\in[1,2]}$ with continuous paths, such that $\tilde X_1=X_1$ and $\tilde X_2=X_2$.
\end{lemma}
%
%
%

\begin{proof} 
By the  martingale version of Skorohod embedding theorem \cite[p.229]{C9}, there exists a Brownian motion $B$ and  optional time $\tau\geq 0$, such that $X_2-X_1=B_{\tau}$ a.s.
Let $M_t=X_1+B_{t\wedge\tau}$. Then by optional stopping theorem $M_t$ is a martingale. This martingale is bounded and, thus, uniformly integrable. Consequently, there exists $M_{\infty}$ such that $\left\{M_t\right\}_{t\in[0,+\infty]}$ is a martingale and  $M_t$ converges to $M_{\infty}$ almost surely by Doob's theorem \cite[p.285]{C8}.

Let $h:[0,+\infty]\rightarrow[1,2]$ be a nondecreasing continuous function and let $\bar M_t\equiv M_{h(t)}$, $t\in[1,2]$. Then $\left\{\bar M_1,\bar M_2\right\}$ has the same joint distribution as $\left\{X_1, X_2\right\}$ and $\left\{\bar M_t\right\}_{t\in[1,2]}$ is  a continuous time martingale with continuous paths.
\end{proof} 
\begin{lemma}
\label{lemma_10} 
Let $\left\{X_t\right\}_{t\in[1,2]}$ be a continuous time martingale with values in $\mathbb R^d$. Let $h:\mathbb R^d\rightarrow\mathbb R$, be a continuous piecewise linear function. Let $\mathbb R^d=\cup_i D_i$, and $h$ is affine on each $D_i$, $i=1,2,\ldots,n$. Then for each $i$
$$E\left[h(X_1)|\left\{X_t\right\}_{t\in[1,2]}\in D_i\right]= E\left[h(X_2)|\left\{X_t\right\}_{t\in[1,2]}\in D_i\right].$$
\end{lemma}
\begin{proof}
The result  immediately follows from conditional Jensen's inequality applied to $h$ and $-h$ restricted to $D_i$ \cite[p.277]{C8}. 
\end{proof}

\noindent Now let us return to the proof of Theorem \ref{th_1}:
\begin{proof}
We start at time $t^*$ and proceed  recursively in the following way.
Let ${B_t}$, $t\in[0,+\infty)$, be a Brownian motion starting at 0. Let ${{\dot X_t}}={X_{t^*}}+{B_{t-t^*}}$ for all $t\geq t^*$. Then ${X_1},{X_2},\ldots,{X_{t^*}},{{\dot X}_t}$ is a martingale. If $\tau_{1}$ is the first hitting time of ${{\dot X_t}}$ on ${\cal K}$, then by the optional stopping theorem ${X_1},{X_2},\ldots,$ ${X}_{t^*},{{\dot X}}_{(t^*+1)\wedge \tau_{1}}$, ${{\dot X}}_{(t^*+2)\wedge \tau_{1}}, \ldots$ is also a martingale. Since this martingale is bounded, it is uniformly integrable, and thus \cite[p.283, p.285]{C8} there exists ${{\dot X}}_{\infty}$, such that $E\left[{{\dot X}}_{\infty}|{\cal B}_n\right]={{\dot X}}_n$ for all $n=t^*,t^*+1,\ldots$, and ${{\dot X}}_n$ converges to ${{\dot X}}_{\infty}$ a.s. Let ${{\tilde X_{t^*}}}={{\dot X}}_{\infty}$. Then ${X_1},{X_2},\ldots,{X_{t^*-1}},{{\tilde X_{t^*}}}$ is a martingale and $P\left({{\tilde X_{t^*}}}\in {\cal K}\right)=1$.

Now let  us make the discrete time martingale ${X_1},{X_2},\ldots,{X_{t^*-1}},{{\tilde X_{t^*}}}$ be continuous on $[t^*-1,t^*]$ and have continuous paths. That is, by Lemma \ref{lemma_9} we can define the martingale $\left\{{{\dot X_t}}\right\}_{t\in[t^*-1,t^*]}$, such that ${{\dot X_{t^*-1}}}= {{X_{t^*-1}}}$, ${{\dot X_{t^*}}}={{\tilde X_{t^*}}}$ and it is continuous. 
Let $\tau_2$ be the stopping time which stops the martingale when it first hits ${\cal K}$ after $t^*-1$. 
Then $\tau_2$ is bounded and $X_{t^*-1},X_{\tau_2}$ is a martingale by the optional stopping theorem.
By Lemma \ref{lemma_10}, ${X_{t^*-1}}\left(={{\dot X_{t^*-1}}}\right)$ and ${{\dot X_{\tau_2}}}$ 
takes the same values on functionals but $P\left({{\dot X_{\tau_{{t^*}-1}}}}\in {\cal K}\right)=1$. 
Since ${X_1},{X_2},\ldots,{X_{t^*-1}},{{\dot X_{\tau_2}}},{{\tilde X_{t^*}}}$ is a martingale, so is  ${X_1},{X_2},\ldots,{X_{t^*-2}},{{\dot X_{\tau_2}}},
{{\tilde X_{t^*}}}$. We let ${{\tilde X_{t^*-1}}}= {{\dot X_{\tau_2}}}$
and proceed in this way until time $1$.
\end{proof}

\noindent We next introduce Problem (b') that is equivalent to Problem
(b).
\\[2pt]

\noindent \parbox[c]{13.5cm}{
\noindent {\bf  Problem (b')}
\hrule
\noindent \rule{0pt}{5mm}
\noindent 
Find the domain of values of
\be
E\left[g(X_{t^*})\right]\nonumber
\ee
under the condition, that $\{X_1, \ldots,X_{t^*}\}$
is a nonnegative martingale with support in ${\cal K}:={\cal K}_{\cal S}\cup{\cal K}_g$ such that
\be
\lb{new_cond2}
E\left[\left(X_t-k_{tj}\right)^+\right] = C_{tj},\quad t\leq t^*,\; j=1,2\ldots, U(t),\\
E\left[\left(X_{t^*}-k_{fut,j}\right)^+\right] \leq C_{fut,j},\quad j=1,2,\ldots, U(_{fut}).\nonumber
\ee
\hrule
\noindent \rule{0pt}{5mm}
}

\noindent From Lemma \ref{l_1} and Theorem \ref{th_1} we obtain 

\begin{corollary}
\label{cor_1}
If the price of the stock is bounded by $L$, as in\rf{bound_on_stock}, then Problem (b) and Problem (b') have the same solutions.
\end{corollary}

\begin{proposition}
If there is no bound on the price of the stock, then as $L\rightarrow\infty$, the optimal solution of Problem (b') converges monotonically to the optimal solution of Problem (b).
\end{proposition}
\begin{proof}
Notice, that since $g$ is a continuous piecewise-linear function, $E\left[g(X_{t^*})\right]$ is always finite once $E[X_{t^*}]$ is finite. Thus, the best possible upper bound on $E\left[g(X_{t^*})\right]$, $C^*$, for Problem (b) is finite . If $C^*_L$ is the best possible upper bound on $E\left[g(X_{t^*})\right]$  for Problem (b'), then $C^*_L\leq C^*$ and $C^*_L$ approaches monotonically  $C^*$ as $L\rightarrow\infty$ according to Corollary \ref{cor_1}. Thus $C^*_L\rightarrow C^*$ as $L\rightarrow\infty$. The same is true for the  lower bound. 
\end{proof}

\subsection{Algorithm for the one dimensional case}

 We have shown that in order to solve problem (b), it is enough to consider Markov martingales with support in ${\cal K}$. The distribution of a discrete time Markov process is fully determined by joint distributions of  each two consecutive states. 

Let ${\cal K}=\left\{k_1,k_2,\ldots,k_{\mathfrak n},\ldots,k_N\right\}$ and ${\EuScript N}$ denote the set of all the indexes $1,2,\ldots,N$.
Then the distribution of the Markov process with support in ${\cal K}$ is determined by the collection of
$$P_{t; {\mathfrak n}_1, {\mathfrak n}_2}=P\left(X_t=k_{{\mathfrak n}_1},X_{t+1}=k_{{\mathfrak n}_2}\right),\;\;t=1,\ldots,t^*,\;\; {\mathfrak n}_1,{\mathfrak n}_2 \in {\EuScript N}.$$

Since most of the constraints that we consider are constraints on marginal distributions of the process, for convenience we also introduce
$$P_{t;{\mathfrak n}}= P\left(X_t=k_{\mathfrak n}\right)\;\;t=1,\ldots,t^*,\;\; {\mathfrak n}\in {\EuScript N}.$$
For $P$ to be a measure, it must satisfy the countable additivity condition:
%
$$\sum_{{\mathfrak n}\in{\EuScript N}} P_{t;{\mathfrak n}_1,{\mathfrak n}}=P_{t;{\mathfrak n}_1}\quad\mbox{for all}\; t\leq t^*,\; {\mathfrak n}_1\in{\EuScript N},$$ 
$$\sum_{{\mathfrak n}\in{\EuScript N}} P_{t;{\mathfrak n},{\mathfrak n}_1}=P_{t+1;{\mathfrak n}_1}\quad\mbox{for all}\; t<t^*,\; {\mathfrak n}_1\in{\EuScript N},$$ 
%
and nonnegativity condition:
$$P_{t;{\mathfrak n}}\geq 0, \;\;\;P_{t;{\mathfrak n}_1,{\mathfrak n}_2}\geq 0, \quad\mbox{for all}\; t<t^*,\;{\mathfrak n},{\mathfrak n}_1,{\mathfrak n}_2\in{\EuScript N}.$$
In addition, for the probability measure we must have
$$\sum_{{\mathfrak n}\in{\EuScript N}} P_{1;{\mathfrak n}}=1.$$
The martingale condition can be written as
$$
 \sum_{{\mathfrak n}\in{\EuScript N}}P_{t;{\mathfrak n}_1,{\mathfrak n}}k_{\mathfrak n}=P_{t;{\mathfrak n}_1}k_{{\mathfrak n}_1}\quad\mbox{for all}\; t\leq t^*,\; {\mathfrak n}_1\in{\EuScript N}.$$
For the martingale to satisfy prices of options with maturities less than or equal to $t^*$, we must have:
$$\sum_{{\mathfrak n}\in{\EuScript N}}P_{t;{\mathfrak n}}(k_{\mathfrak n}-k_{tj})^+=C_{tj}
\quad\mbox{for all}\;t\leq t^*,\; j=1,2\ldots, U(t)$$
As was proved in Theorem \ref{th_1}, for prices of options of maturities greater than $t^*$ to be satisfied, it is enough to have
$$\sum_{{\mathfrak n}\in{\EuScript N}}P_{t^*;{\mathfrak n}}(k_{\mathfrak n}-k_{fut,j})^+\leq C_{fut,j}
\quad\mbox{for all}\;j=1,2\ldots, U(_{fut}).$$

Thus the problem of determining best possible upper (lower) bounds on the price of an option with maturity $t^*$ and a payoff function $g$, consists of finding the collection of $P_{t;{\mathfrak n}}$ and $P_{t;{\mathfrak n}_1,{\mathfrak n}_2}$, such that 
$$
\begin{array}{rl}
\max (\min) & \displaystyle 
\sum_{{{\mathfrak n}\in {\EuScript N}}}P_{t^*;{\mathfrak n}}\;
g(k_{\mathfrak n}) \vspace{3pt}\\
{\rm s.t.} & \displaystyle \sum_{{{\mathfrak n}\in {\EuScript N}}} P_{1;{\mathfrak n}}=1, \vspace{3pt}\\
&\displaystyle 
\sum_{{\mathfrak n}\in{\EuScript N}} P_{t;{\mathfrak n}_1,{\mathfrak n}}=P_{t;{\mathfrak n}_1}\quad\mbox{for all}\; t\leq t^*,\; {\mathfrak n}_1\in{\EuScript N}, \vspace{3pt}\\
&\displaystyle \sum_{{\mathfrak n}\in{\EuScript N}} P_{t;{\mathfrak n},{\mathfrak n}_1}=P_{t+1;{\mathfrak n}_1}\quad\mbox{for all}\; t<t^*,\; {\mathfrak n}_1\in{\EuScript N}, \vspace{3pt}\\
&\displaystyle \sum_{{\mathfrak n}\in{\EuScript N}}P_{t;{\mathfrak n}_1,{\mathfrak n}}k_{\mathfrak n}=P_{t;{\mathfrak n}_1}k_{{\mathfrak n}_1}\quad\mbox{for all}\; t<t^*,\; {\mathfrak n}_1\in{\EuScript N}, \vspace{3pt}\\
&\displaystyle \sum_{{\mathfrak n}\in{\EuScript N}}P_{t;{\mathfrak n}}(k_{\mathfrak n}-k_{tj})^+=C_{tj}\;
\quad\mbox{for all}\;t\leq t^*,\; j=1,2\ldots, U(t), \vspace{3pt}\\
&\displaystyle \sum_{{\mathfrak n}\in{\EuScript N}}P_{t^*;{\mathfrak n}}(k_{\mathfrak n}-k_{fut,j})^+\leq C_{fut,j}\quad\mbox{for all}\;j=1,2\ldots, U(_{fut}),\vspace{3pt}\\
& P_{t;{\mathfrak n}_1,{\mathfrak n}_2}\geq 0\quad\mbox{for all}\; t\leq t^*,\; {\mathfrak n}_1,{\mathfrak n}_2\in{\EuScript N}.
\end{array}
$$
This is a linear optimization problem with  $O(N^2)$, where $N$ is the number of elements of ${\cal K}$.

\section{The Two Dimensional Case}
In this section, we address Problem (c), i.e., 
 the problem of determining bounds on options on two different assets if prices of options on individual assets are known.

Let $\left\{S^{^I}(t)\right\}_{t\geq 0}$ be a stochastic process describing the price process of the first asset, and $\left\{S^{^{II}}(t)\right\}_{t\geq 0}$ that of the second asset.
If $\left(k_{ij}^{{h}},C_{ij}^{{h}},t_i\right)$, $i=1,\ldots,n$, $j=1,2,\ldots, U(t,{h})$,  represent the set of options on asset $S^{{h}}$, with $t_i$ being the maturity of an option, $k_{ij}^{{h}}$ its strike price, and $C_{ij}^{{h}}$ its price, then we must have 
\be
\lb{first}
C_{ij}^{{h}}=e^{\displaystyle{-\int_0^{t_i}r(s)ds}}E_{Q}\left[\left({S^{{h}}({t_i})-k_{ij}^{{h}}}\right)^+\right]=E_{Q}\left[\left({X_{t_i}^{{h}}-k_{ij}^{{h}}e^{\displaystyle{-\int_0^{t_i}r(s)ds}}}\right)^+\right], 
\ee
where $\left\{X_{t}^{{h}}\right\}_{t\geq 0}$, is a martingale under $Q$. 
As in the single dimensional case, we can assume without loss of generality that 
 the risk-free rate $r(s)$ can be assumed 0. Thus, we can change times $t_i$ to their indexes.

Let ${\cal
T}:=\left\{({h},t,j)|{h}\in\{_{^I},_{^{II}}\},\;t\in\{1,2,\ldots,n\},
j=1,2,\ldots U(t,{h})\right\}$. We define a family of measurable functions $f_{tj}^{{h}}: \mathbb R^{2}_+\rightarrow \mathbb R_+$ as
\be
f_{tj}^{{h}}({\bold X_t})=\left(X_t^{{h}}-k_{tj}^{{h}}\right)^+,\;({h},t,j)\in {\cal T}.
\ee
Then the problem can be formulated as follows:
\\[2pt]

\noindent \parbox[c]{13.5cm}{
\noindent {\bf  Problem (c)}

\noindent
\hrule
\noindent \rule{0pt}{5mm}
\noindent 
Given the set $\left(k_{tj}^{{h}},C_{tj}^{{h}}\right)$, $({h},t,j)\in{\cal T}$, and a continuous piecewise linear function $g:\mathbb R^{2}\rightarrow \mathbb R$, 
find the maximal and minimal possible values of  $E[g({\bold X_{t^*}})]$, for a given $t^*\in\{1,2,\ldots,n\}$, under the condition that  ${\bold X_1},{\bold X_2},\ldots,{\bold X_n}$ is a two dimensional nonnegative martingale,  such that
\be
\lb{conditions}
E\left[f_{tj}^{{h}}({\bold X_t})\right]=C_{tj}^{{h}},\;({h},t,j)\in {\cal T}
\ee
\hrule
\noindent \rule{0pt}{5mm}
}

\begin{definition}
\label{def_9_2}
We call function $g:\mathbb R^d\rightarrow \mathbb R$ piecewise-linear, if there exists a partition of $\mathbb R^d$ by a finite number of hyperplanes (lines if $d=2$) into nonoverlapping subsets (regions if $d=2$), such that in the interior of each subset $g$ is an affine function. 
\end{definition}

Notice, that by Proposition \ref{prop_6_2} we might restrict our consideration to Markov martingales only.

%
\subsection{Treatment of future conditions}
%
Notice that the function $f_{tj}^{{h}}$ depends only on the ${h}$-coordinate of ${\bold X_t}$. With this property of $f-$constraints, we can show that future conditions have a ``simple effect'',
in the sense that conditions\rf{conditions} for $t>t^*$ can be reduced to separate conditions on the first and second coordinates of ${\bold X_{t^*}}$.
Let
$${\cal S}_{\geq t^*}^{{h}}=\left\{\left.\left(k_{tj}^{{h}},C_{tj}^{{h}}\right)\right|t^*\leq t\leq n, t\in\{1,2,\ldots,U(t,{{h}})\} \right\},\;{h}=_{^I},_{^{II}},$$
$${\cal S}_{\geq t^*}^{h\infty} = {\cal S}_{\geq t^*}^{{h}}\cup \{{{(0,+\infty)},(+\infty,0)}\},\;{h}=_{^I},_{^{II}},$$
and
\be
{\cal F}^{{h}}=\left\{\left(k_{tj}^{{h}},C_{tj}^{{h}}\right)|t^*<t\leq n,\;\left(k_{tj}^{{h}},C_{tj}^{{h}}\right)\;
{\mbox{is a vertex of the convex hull of }}{\cal S}_{\geq t^*}^{h\infty}\right\}\nonumber
\ee
 represent future points which are on the border of the convex hull of ${\cal S}_{\geq t^*}^{h\infty}$. 
Let 
$$\left(k_{fut,j}^{{h}},C_{fut,j}^{{h}}\right),\quad  j=1,2,\ldots,U(_{fut},{h})$$ 
be enumerated elements of ${\cal F}^{{h}}$.

In formulating all further results we will assume that $\left(k_{tj}^{^I},C_{tj}^{^I}\right)$,\linebreak
 $(_{^I},t,j)\in {\cal T}$ and $\left(k_{tj}^{^{II}},C_{tj}^{^{II}}\right)$, $(_{^{II}},t,j)\in {\cal T}$ separately do satisfy the no-arbitrage condition, that is that there exist one-dimensional martingales $X_t^{^I}$ and $X_t^{^{II}}$ such that
$$E\left[\left(X_t^{{h}}-k_{tj}^{{h}}\right)^+\right]=C_{tj}^{{h}}, \quad{h}=_{^I},_{^{II}}.$$
Notice that this condition ensures that the two dimensional martingale which satisfies conditions\rf{conditions} also exists.

\begin{theorem}
\label{th6}
 Let ${\bold X_1},{\bold X_2},\ldots,{\bold X_{t^*}}$ be a two-dimensional martingale satisfying\\ 
$E\left[f_{tj}^{{h}}({\bold X_t})\right]=C_{tj}^{{h}}$, for all $\{({h},t, j)\in{\cal T}\;|\;t\leq t^*\}$. If also for each  ${h}\in\{_{^I},_{^{II}}\}$ 
$$E\left[\left(X^{{h}}_{t^*}-k_{fut,j}^{{h}}\right)^+\right]\leq C_{fut,j}^{{h}}$$ is satisfied for all $(k_{fut,j}^{{h}},C_{fut,j}^{{h}})\in {\cal F}^{{h}}$, then there exist two-dimensional random variables ${\bold X_{t^*+1}},\ldots,{\bold X_n}$ such that 
\begin{itemize}
\item[{\bf (a)}] ${\bold X_1},{\bold X_2},\ldots,{\bold X_{t^*}},{\bold X_{t^*+1}},\ldots,{\bold X_n}$ is a martingale;
\item[{\bf (b)}] all the conditions $E\left[f_{tj}^{{h}}({\bold X_t})\right]=C_{tj}^{{h}}$ are satisfied for all $({h},t,j)\in {\cal T}$. 
\end{itemize}
\end{theorem}
%
%
\begin{proof}
As in the proof of Lemma \ref{l_1}, we can show that for each  ${h}\in\{_{^I},_{^{II}}\}$ there exists a one-dimensional martingale 
  $\tilde X_{t^*}^{{h}},\tilde X_{t^*+1}^{{h}},\ldots,\tilde X_n^{{h}}$, such that it satisfies the set $\left(k_{tj}^{{h}},C_{tj}^{{h}}\right)$, $t\geq t^*,\; j=1,\ldots,U(t,{h})$ and $\Psi_{\tilde X_{t^*}^{{h}}}=\Psi_{X^{{h}}_{t^*}}$.

%
%
 

We will extend the martingale ${\bold X_1},{\bold X_2},\ldots,{\bold X_{t^*}}$ to the martingale ${\bold X_1},{\bold X_2},\ldots,{\bold X_{t^*}},$ ${\bold X_{t^*+1}}, \ldots,{\bold X_n}$, such that ${\bold X_{t^*}},{\bold X_{t^*+1}},\ldots,{\bold X_n}$ is a Markov process. To define a joint distribution of ${\bold X_{t}}$ and ${\bold X_{t+1}}$ for $t\geq t^*$, allow each coordinate to evolve independently of each other. That is, for each ${{h}}=_{^I},_{^{II}}$ define ${\cal L}\left(X^{{h}}_{t}|X^{{h}}_{t+1}\right)$, as ${\cal L}\left(\tilde X^{{h}}_{t}|\tilde X^{{h}}_{t+1}\right)$ for $t=t^*,t^*+1,\ldots,n$. 
\end{proof}

Theorem \ref{th6} allows us to reduce ``future conditions'' to conditions at time $t^*$. This is a significant simplification of the problem, since difficulties arise from  the requirement that the martingale satisfies $E\left[f_{tj}^{{h}}({\bold X_t})\right]=C_{tj}^{{h}},\;({h},t,j)\in {\cal T}$.
\begin{corollary}
\label{cor_5}
Our main problem is equivalent to the following:
find the maximum (minimum) of $E\left[g\left({\bold X_{t^*}}\right)\right]$ under the condition that ${\bold X_1},{\bold X_2},\ldots,{\bold X_{t^*}}$ is a Markov Martingale such that 
\\[1pt]
{\bf (a)} $E\left[\left(X_t^{{h}}-k_{tj}^{{h}}\right)^+\right]=C_{tj}^{{h}}$, for all $t\leq t^*$, $j\leq U(t,{h})$ and
\\[1pt]
{\bf (b)} $E\left[\left(X_{t^*}^{{h}}-k_{fut,j}^{{h}}\right)^+\right]\leq C_{fut,j}^{{h}}$ for all $\left(k_{fut,j}^{{h}},C_{fut,j}^{{h}}\right)\in{\cal F}^{{h}}$. 
\end{corollary}
\begin{definition}
\label{def_11}
Let $f^{{h}}_{fut,j}:\mathbb R^{2}_+\rightarrow\mathbb R_+$, $f^{{h}}_{fut,j}({\bold X})=\left(X^{{h}}-k_{fut,j}^{{h}}\right)^+$, be a family of functions, where each function corresponds to some $k_{fut,j}^{{h}}$, such that  $\left(k_{fut,j}^{{h}},C_{fut,j}^{{h}}\right)\in {\cal F}^{{h}}$, ${h}=_{^I},_{^{II}}$. 
\end{definition}
Constraints from the past introduce much more difficulty to the problem than future constraints. It is not possible to reduce  past martingale conditions to separate  conditions on $X^{^I}_{t^*}$ and $X^{^{II}}_{t^*}$. 
In other words, for a given  distribution of ${\bold X_{t^*}}$, and two sequences of marginal laws $\left\{\Pi_t\right\}_{t=1,\ldots,t^*}$ and $\left\{Q_t\right\}_{t=1,\ldots,t^*}$  of one dimensional processes, such that $\Pi_{t^*}={\cal L}\left(X^{^{I}}_{t^*}\right)$ and $Q_{t^*}={\cal L}\left(X^{^{II}}_{t^*}\right)$, the existence of one-dimensional martingales with the above sequences of marginal laws, is not sufficient for the existence of a two-dimensional martingale ${\bold X_1},{\bold X_2},\ldots,{\bold X_{t^*}}$ such that ${\cal L}\left(X^{^{I}}_{t}\right)=\Pi_t$ and  ${\cal L}\left(X^{^{II}}_{t}\right)=Q_t$ for $t=1,\ldots,t^*$. 

\begin{figure}[htbp]
\begin{center}
\input{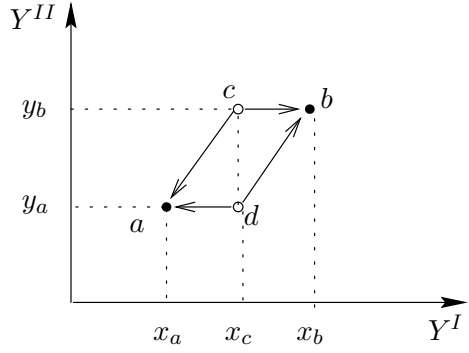}
\caption{Counterexample.}
\label{counterexample}
\end{center}
\end{figure}

\noindent
{\bf Counterexample}. Suppose that ${\bold X_{2}}$ has an atomic distributions with two atoms $a$ and $b$, as shown on Figure \ref{counterexample}, where each atom has weight of 1/2. Suppose that  the marginal law of $X^{^I}_1$ is $P\left(X^{^I}_1=x_c\right)=1$, where $x_c=\frac12(x_a+x_b)$ and the marginal law of $X^{^{II}}_1$ is  $P\left(X^{^{II}}_1=y_a\right)=P\left(X^{^{II}}_1=y_b\right)=\frac12$. Obviously, we can define a joint distribution of  $X^{^{I}}_1, X^{^{I}}_2$, so that it is a one-dimensional martingale. The same is true for $X^{^{II}}_1, X^{^{II}}_2$. However, it is not possible to define a joint distribution of ${\bold X_1}$, ${\bold X_2}$, so that it is a two-dimensional martingale.

Thus, the condition that ${\bold X_1},{\bold X_2},\ldots,{\bold X_{t^*}}$ is a martingale, is truly a ``joint'' condition on $\left\{X^{^{I}}_t\right\}$ and $\left\{X^{^{II}}_t\right\}$. 
\subsection{Restriction to discrete distributions}
In what follows, function $g$ is defined as
$g:\mathbb R^{2}_+\rightarrow \mathbb R_+,\;(x,y)\rightarrow\left|\alpha x+\beta y-k_g\right|^+$. However, an analogous solution would apply if $g:\mathbb R^2_+\rightarrow \mathbb R$ is any continuous piecewise-linear function.

As in the one dimensional case we add an additional constraint  that stock prices are bounded, that is, there exists $L>0$, such that for ${h}=I,II$
$$P\left(X_1^{h},X_2^{{h}},\ldots, X_n^{{h}}\leq L\right)=1.$$

\begin{definition} 
\label{Definition_Graph_G}
For each $t=1,2,\ldots,t^*-1$ we define Graph $G_t$ 
(Figure \ref{graph_G_1}), which is contained in $[0,L]\times[0,L]$ 
and has vertices on intersections of the lines:
\\[1pt]
{\bf (a)} $x=k_{tj}^{^I}$, $j\leq U(t,_{^I})$,
$\quad$ {\bf (b)} $y=k_{tj}^{^{II}}$,  $j\leq U(t,_{^{II}})$,
\\[1pt]
{\bf (c)} segments of $Q$: from $(0,0)$ to $(0,L)$, $(0,L)$ to $(L,L)$, $(L,L)$ to $(L,0)$, $(0,0)$ to $(L,0)$.
\\[1pt]

\noindent For $t^*$ we define Graph $G_{t^*}$ (Figure \ref{graph_G_2}) which has vertices on intersections of the lines
\\[1pt]
{\bf (a)} $x=k_{t^*j}^{^I}$,  $j\leq U(t,_{^I})$,$\;\;\;\;\;$ ($\mbox{{\bf a}}^{\prime}$) $x=k^{^I}_{fut,j}$, $j\leq U(_{fut},_{^I})$,
\\[1pt]
{\bf (b)} $y=k_{t^*j}^{^{II}}$,  $j\leq U(t,_{^{II}})$,$\;\;\;$ ($\mbox{{\bf b}}^{\prime}$) $y=k^{^{II}}_{fut,j}$, $j\leq U(_{fut},_{^{II}})$,
\\[1pt]
{\bf (c)} segments of $Q$: from $(0,0)$ to $(0,L)$, $(0,L)$ to $(L,L)$, $(L,L)$ to $(L,0)$, $(0,0)$ to $(L,0)$,
\\[1pt]
{\bf (d)} $\alpha x+\beta y=k_g$.

\end{definition}
Each edge of $G_t$  belongs to one of the corresponding  lines (or segments). The intersection between two edges is either empty or consists of one vertex point.

\begin{figure}[h]
  \begin{minipage}[t]{.45\textwidth}
    \begin{center}  
     \input{apr_1_2.pstex_t}
      \caption{The form of Graph $G_t$ for $t=1,2,\ldots,t^*-1$.}
      \label{graph_G_1}
    \end{center}
  \end{minipage}
  \hfill
  \begin{minipage}[t]{.45\textwidth}
    \begin{center}  
      \input{apr_1_3.pstex_t}
      \caption{The form of Graph $G_{t^*}$.}
      \label{graph_G_2}
    \end{center}
  \end{minipage}
  \hfill
\end{figure}

Let 
\begin{itemize}
\item 
$G_t^V$ designate all the vertices of Graph $G_t$,
\item
$G_t^E$ denote the set of points which lie on edges of $G_t$ and
\item
$G_t^R$  denote regions of $\mathbb R^2$ over which edges of Graph  $G_t$ splits $[0,L]\times[0,L]$.
\end{itemize}

Let us enumerate regions of $G_t^R$. We will thus write $r$ for the region of $G_t^R$ numbered $r$.  Let ${\EuScript R}_t$  be the set of all integers such that each of them represents some region of $G_t^R$ and  ${\EuScript V}$  be the set of all integers such that each of them represents some vertex of $G_{t^*}^V$. 
Let ${\mathfrak {Reg}}_t: [0,L]\times[0,L]\rightarrow {\EuScript R}_t$, be the function which attributes to each point of $[0,L]\times[0,L]$ the index of the  region that contains it, and let  ${\EuScript{V}}{\it {er}}: [0,L]\times[0,L]\rightarrow {\EuScript V}$, be the function which attributes to each vertex of $G_{t^*}^V$ its index $v$.  Since edges are contained by several regions, we have to specify which region corresponds to each edge.

\noindent {\bf Definition 13.} We will say that a random process ${\bold X_1},{\bold X_2},\ldots,{\bold X_{t^*}}$ has property $P_G$ if and only if all of the following hold
\\[1pt]
{\bf (a)} $P\left({\bold X_{t^*}}\in G^V_{t^*}\right)=1$;
\\[1pt]
{\bf (b)} ${\bold X_1}$ has an atomic distribution with at most one atom per region of $G^R_1$;
\\[1pt]
{\bf(c)} If ${r}_1\in {\EuScript R}_1, \;{r}_2\in {\EuScript R}_2\ldots,{r}_t\in{\EuScript R}_t$ with $t<t^*$, then 
$${\cal L}\left({\bold X_{t+1}}\left|\left[{\mathfrak {Reg}}_t({\bold X_1}),{\mathfrak {Reg}}_t({\bold X_2}),\ldots,{\mathfrak {Reg}}_t({\bold X_t})\right.\right]=[r_1,r_2,\ldots,r_t]\right)$$ 
has an atomic measure with at most one atom per region of $G^R_t$.
\begin{theorem}
 \label{th_graph2}
Let ${\bold X_1},{\bold X_2},\ldots,{\bold X_{t^*}}$ be a nonnegative Markov Martingale bounded by $L$. Then there exists a nonnegative Markov Martingale $\left\{{\bold {\bar X_1}},{\bold {\bar X_2}},\ldots,{\bold {\bar X_{t^*}}}\right\}$, bounded by $L$, such that it has property $P_G$ and satisfies

$$E\left[f_{tj}^{{h}}({\bold X_t})\right]=E\left[f_{tj}^{{h}}\left({\bold {\bar X_t}}\right)\right],\; t\leq t^*,\; ({h},t,j)\in {\cal T},$$%
$$E\left[f_{fut,j}^{{h}}({\bold X_{t^*}})\right]=E\left[f_{fut,j}^{{h}}\left({\bold {\bar X_{t^*}}}\right)\right],\; j=1,2,\ldots,U(_{fut},{h}),\; {h}\in\{_{^I},_{^{II}}\},$$
$$E\left[g\left({\bold X_{t^*}}\right)\right]=E\left[g\left({\bold{\bar X_{t^*}}}\right)\right].$$ 
\end{theorem}
\begin{proof}
For $i = 1,\ldots,t^*$, let $\sigma_i$ designate  the smallest $\sigma$-algebra for which ${\bold X_1}, {\bold X_2},\ldots,{\bold X_i}$ are all measurable. 

The proof goes by induction. Let $\bar\sigma_1$ be the smallest $\sigma$-algebra for which the random variable ${\mathfrak {Reg}_1}({\bold X_1})$ is measurable.  Then $\bar\sigma_1\subset\sigma_1$. Define ${\bold{\bar X}_1} = E\left[{\bold X_1}|\bar\sigma_1\right]$. Then $\left\{{\bold{\bar X}_1},{\bold X_1};\bar\sigma_1, \sigma_1\right\}$ is a martingale. ${\bold {\bar X}_1}$ has an atomic distribution  with at most one atom ${\bold Z_{1;{r}}}= E\left[{\bold {X_1}}\left|{\mathfrak
{Reg}_1}\left({\bold{X_1}}\right)={r}\right.\right]$ per each region $r$ of $G^R_1$ and takes the same values on functionals $f_{1j}^{{h}}$, $j=1,2,\ldots,U(1,{h})$, ${h}\in\{_{^I},_{^{II}}\}$, as ${\bold {X_1}}$ by Lemma \ref{lemma_10}.

Now, ${\bold {\bar X}_1},{\bold{X_2}},\ldots, {\bold {X_{t^*}}}$ is also a martingale. Let $\bar \sigma_2$ be the smallest $\sigma$-algebra which includes $\bar\sigma_1$ and for which ${\mathfrak {Reg}_2}({\bold X_2})$ is measurable. Define
${\bold{\bar X}_2} = E\left[{\bold X_2}|\bar\sigma_2\right]$.  Then  $\left\{{\bold{\bar X}_1},{\bold {\bar X}_2}, {\bold {X}_2} ;\bar\sigma_1, \bar \sigma _2, \sigma_2\right\}$ is a martingale. Indeed, $E\left[{\bold {\bar X}_2}|\bar\sigma_1\right] = E\left[E\left[{\bold{X}_2}|\bar\sigma_2\right]|\bar\sigma_1\right]=$ 

\noindent 
$E\left[{\bold{X}_2}|\bar\sigma_1\right] = {\bold {\bar X}_1}$, since $\bar\sigma_1\subset\bar\sigma_2$, and $E\left[{\bold X_2}|\bar\sigma_2\right]={\bold{\bar X}_2}$ by definition.
${\cal L}\left({\bold {\bar X_2}}\left|{\mathfrak {Reg}}_1\left({\bold{X_1}}\right)={r}\right.\right)$ has an atomic measure with at most one atom 
${\bold Z_{2;{r}, \tilde r}}= E\left[{\bold {X_2}}\left|{\mathfrak
{Reg}_1}\left({\bold{X_1}}\right)={r},{\mathfrak
{Reg}_2}\left({\bold{X_2}}\right)={\tilde r} \right.\right]$
per each region $\tilde r$ of $G_2^R$.
Also by Lemma  \ref{lemma_10} ${\bold{\bar X_2}}$ takes the same values on the functionals  $f_{2j}^{{h}}$, $j=1,2,\ldots,U(1,{h})$, ${h}\in\{_{^I},_{^{II}}\}$,  as ${\bold X_2}$.

Continuing in the same way, 
we can form a Markov Martingale ${\bold{\bar X_1}},\ldots,{\bold{\bar X_{t^*-1}}}, {\bold{\bar X_{t^*}}}$ 
such that  
${\cal L}\left({\bold{\bar X_t}}|\left[{\mathfrak {Reg}}_1({\bold X_1}),
{\mathfrak {Reg}}_2({\bold X_2}),\ldots,{\mathfrak {Reg}}_{t-1}({\bold X_{t-1}})\right]=[{r}_1,{r}_2,\ldots,{r}_{t-1}]\right)$ has an atomic measure with at most one atom per region of ${\EuScript R}_t$ for all $t=1,2,\ldots,t^*$. 

Now, let ${\bold B}_{t-t^*}$ be a two-dimensional Brownian motion starting at 0. 
We define for $t\geq t^*$ ${\bold{\hat X}}_t = {\bold{\bar X_{t^*}}}+{\bold B}_{t-t^*}$ and let $\tau$ be the first hitting time of ${\bold{\hat X}}_t$ on $G^E_{t^*}$. Then ${\bold{\bar X}}_{t^*}, {\bold{\hat X}}_{t^*\wedge\tau}, {\bold{\hat X}}_{(t^*+1)\wedge\tau},\ldots$, is a martingale by the optional stopping theorem. This martingale is bounded, thus there exists ${\bold{\hat X}}_{\infty}$, such that ${\bold{\hat X}}_{(t^*+1)\wedge\tau}, {\bold{\hat X}}_{(t^*+2)\wedge\tau},\ldots,{\bold{\hat X}}_{\infty}$ is a martingale and  ${\bold{\hat X}}_{(t^*+1)\wedge\tau}, {\bold{\hat X}}_{(t^*+2)\wedge\tau},\ldots,$ converges to ${\bold{\hat X}}_{\infty}$ a.s. by Doob's theorem. Let ${\bold{\dot X}}_{t^*}={\bold{\hat X}}_{\infty}$.

By Lemma \ref{lemma_10}, $P\left({\bold{\dot X}}_{t^*}\in G^E_{t^*}\right)=1$ and
${\bold{\dot X}}_{t^*}$ takes the same values on functionals $E[f_{t^*j}^{{h}}(.)]$, $E[f_{fut,j}^{{h}}(.)]$ and $E[g(.)]$ as ${\bold{\bar X_{t^*}}}$.
In the same way, but by adding an one dimensional Brownian motion to ${\bold{\dot X}}_{t^*}$ over the edge to which ${\bold{\dot X}}_{t^*}$ belongs, we can reduce the support of the martingale at time $t^*$ to $G^V_{t^*}$.
\end{proof}

Thus, to solve  optimization problem (c), it is enough to consider only Markov martingales ${\bold X_1},{\bold X_2},\ldots,{\bold X_{t^*}}$ satisfying the property $P_G$. Consequently, it is a finite dimensional optimization problem (distributions with the $P_G$ property can be described with finitely many parameters.) 
General measures belong to an infinite-dimensional space which is even uncountable. So the above result is a significant simplification of the 
problem. However we will see that the number of parameters grows exponentially with $t^*$. So solution methods described below  will work only if $t^*$ is not
 very large.

\subsection{Algorithm for the two dimensional case}

By Theorem \ref{th_graph2} in order to solve Problem (c) it is enough to consider martingales with property $P_G$. That means that the distribution of the initial state ($t=1$) is atomic with at most one atom per each region of 
 $G_1^R$, the distribution of the final state ($t=t^*$) is atomic with atoms at the vertices of  $G_{t^*}$, and the conditional distribution of each intermediate state ($1<t<t^*$) is atomic with at most one atom per each region of $G_t^R$, conditioned that the sequence of all prior visited regions is known.

\begin{figure}[h]
\begin{center}
\input{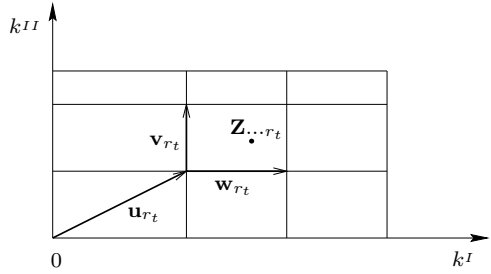}
\caption{Parametrization of a rectangular.}
\label{rect_1}
\end{center}
\end{figure}

The distribution of the described process is fully determined by the probabilities
\be
P_{t^*;r_1,\ldots,r_{t^*-1},r_{t^*}} = P\left(\left[{\mathfrak {Reg}}_1({\bold X_1}),\ldots,{\mathfrak {Reg}}_{t^*-1}({\bold X_{t^*-1}}), {\EuScript V}{\it {er}}({\bold X_{t^*}})\right]=[r_1,\ldots,r_{t^*-1}, v]\right)
\ee
of visiting all possible sequences of regions (vertices at the final state), and  by the unique places of visits of each particular region, conditioned on the prior sequence of regions. 
Let ${\bold  Z_{t;r_1,r_2,\ldots,r_t}}\in \mathbb R^{2}_+$ denote the place of the unique atom  of the conditional distribution 
\be
{\cal L}\left({\bold X_t}\;\left|\;\left[{\mathfrak {Reg}}_1({\bold X_1}),\ldots,{\mathfrak {Reg}}_t({\bold X_t})\right]=[r_1,r_2,\ldots,r_t]\right.\right),\nonumber
\ee
where  $r_1\in{\EuScript R}_1,\;  r_2\in{\EuScript R}_2,\ldots,  r_t\in{\EuScript R}_t$, $t<t^*$.
Since all regions of ${\EuScript R}_t$ for $t=1,2,\ldots,t^*-1$ are rectangular,
${\bold  Z_{t;r_1,r_2,\ldots,r_t}}$ can be represented as 
\be
\lb{parametrization}
{\bold  Z_{t;r_1,r_2,\ldots,r_t}}={\bold u}_{r_t}+{\bold w}_{r_t}z^{_I}_{t;r_1,r_2,\ldots,r_t}+{\bold v}_{r_t}z^{_{II}}_{t;r_1,r_2,\ldots,r_t},
\ee
where ${\bold u}_{r_t},{\bold w}_{r_t},{\bold v}_{r_t}\in\mathbb R^{2}_+$ are vectors uniquely defined by the rectangular $r_t$ (see Figure \ref{rect_1}) and 
$z^{{h}}_{t;r_1,r_2,\ldots,r_t}\in[0,1]$, ${h}=_{^I},_{^{II}}$. 

Notice that $f_{tj}^{^I}$, $f_{tj}^{^{II}}$ are affine when restricted to one of the regions of $G_t^R$, and 
\be
\lb{function_decomposition}
f_{tj}^{{h}}\left({\bold Z}_{t;r_1,\ldots,r_{t-1},r}\right)={a_{tj}^{{h}}}(r) z^{{h}}_{t;r_1,r_2,\ldots,r_{t-1},r}+ b_{tj}^{{h}}(r),\quad r\in{\EuScript R}_t
\ee
for some coefficients $a_{tj}^{{h}}(r)$, $b_{tj}^{{h}}(r)$, 
which are easy to determine.

For convenience we also introduce 
$$P_{t;r_1,\ldots,r_{t}} = P\left(\left[{\mathfrak {Reg}}_1({\bold X_1}),\ldots,{\mathfrak {Reg}}_{t}({\bold X_{t}})\right]=[r_1,\ldots,r_{t}]\right)
$$
for each $t<t^*$, which are probabilities of visiting a particular sequence of regions prior to time $t$.

In order to satisfy the conditions  of a probability measure, in addition to nonnegativity, we must have
\be
\lb{prob_meas_cond}
\sum_{r\in{\EuScript R}_1} P_{1;r}=1,
\ee

\be
P_{t-1;r_1,r_2,\ldots,r_{t-1}}=\sum_{r\in{\EuScript R}_t} P_{t;r_1,r_2,\ldots,r_{t-1},r},\;\mbox{for all}\;\;r_1\in{\EuScript R}_1,\ldots,r_{t-1}\in{\EuScript R}_{t-1}, \;t=2,\ldots,t^*-1,
\ee
\be
P_{t^*-1;r_1,r_2,\ldots,r_{t^*-1}}=\sum_{r\in{\EuScript R}_t} P_{t;r_1,r_2,\ldots,r_{t^*-1},v},\;\; \mbox{for all}\;\;r_1\in{\EuScript R}_1,\ldots,r_{t^*-1}\in{\EuScript R}_{t^*-1},\;v\in {\EuScript V}.
\ee

\noindent
To satisfy the martingale condition we must have
\be
&{\bold Z_{t-1;r_1r_2,\ldots,r_{t-1}}}P_{t-1;r_1,r_2,\ldots,r_{t-1}}=\sum_{r\in{\EuScript R}}{\bold Z_{t;r_1r_2,\ldots,r_{t-1},r}}P_{t;r_1,r_2,\ldots,r_{t-1},r}\\
& \mbox{for all}\;\;r_1\in{\EuScript R}_1,\ldots,r_{t-1}\in{\EuScript R}_{t-1}, \;\; t=2,3,\ldots,t^*-1,\nonumber
%
\ee
\be
&{\bold Z_{t^*-1;r_1r_2,\ldots,r_{t^*-1}}}P_{t-1;r_1,r_2,\ldots,r_{t^*-1}}=\sum_{v\in{\EuScript V}}{\bold Z_{t^*;r_1r_2,\ldots,r_{t^*-1},v}}P_{t^*;r_1,r_2,\ldots,r_{t-1},v}\\
& \mbox{for all}\;\;r_1\in{\EuScript R}_1,\ldots,r_{t^*-1}\in{\EuScript R}_{t^*-1},\;v\in {\EuScript V}.\nonumber
%
\ee
Notice, that $\bold Z_{t;r_1r_2,\ldots,r_{t-1},r}$ is defined in terms of $z^{_I}_{t;r_1r_2,\ldots,r_{t-1},r}$ and $z^{_{II}}_{t;r_1r_2,\ldots,r_{t-1},r}$ in\rf{parametrization}.
\\[5pt]
Finally, the measure must satisfy the prices of options. We use the representation \rf{function_decomposition} of payoff functions to write down the constraints from options prices. For satisfying past constraints we must have 
\be
&\sum_{r_1\in{\EuScript R}_1,\ldots,r_{t}\in{\EuScript R}_{t}}\left(a_{tj}^{{h}}(r_t) z^{{h}}_{t;r_1,r_2,\ldots,r_{t-1},r_t}+{b_{tj}}^{{h}}(r_t)\right)
P_{t;r_1,r_2,\ldots,r_t}=C_{tj}^{{h}},\\
& {\mbox {for }}\left\{({h},t,j)\in {\cal T}|\;t<t^*.\right\}\nonumber
\ee
The present constraint are
\be
&\sum_{r_1\in{\EuScript R}_1,\ldots,r_{t^*-1}\in{\EuScript R}_{t^*-1},v\in  {\EuScript V}}f_{t^*j}^{{h}}(v)
P_{t^*;r_1,\ldots,r_{t^*-1},v}=C_{t^*j}^{{h}}, \;\;\\
& {\mbox {for }}\left\{({h},t,j)\in {\cal T}|\;t=t^*.\right\}\nonumber
\ee
Finally, constraints from options of maturities larger than $t^*$ are
\be
\label{fc1}
&\sum_{r_1\in{\EuScript R}_1,\ldots,r_{t^*-1}\in{\EuScript R}_{t^*-1},v\in  {\EuScript V}}f_{fut,j}^{{h}}(v)
P_{t^*;r_1,r_2,\ldots,r_{t^*-1},v}\leq C_{fut,j}^{{h}}\\
& {\mbox {for }}\left\{({h},t,j)\in {\cal T}|\;t>t^*.\right\}
\nonumber
\ee

Thus problem (c) can be written as : 
\\[1pt]
\hrule
\noindent \rule{0pt}{5mm}
\noindent
Find the collection of  $P_{t;r_1,r_2,\ldots,r_t}\geq 0$ and $z^{{h}}_{t;r_1,\ldots,r_t}\in[0,1]$, $r_1\in{\EuScript R}_1,\ldots,r_t\in{\EuScript R}_t,\;t=1,2,\ldots,{t^*}-1$, ${h}= _{^I},_{^{II}}$, such that it satisfies conditions\rf{prob_meas_cond}-(\ref{fc1}) and  
$$\sum_{r_1\in{\EuScript R}_1,\ldots,r_{t^*-1}\in{\EuScript R}_{t^*-1},v\in  {\EuScript V}}g(v)P_{t^*;r_1,\ldots,r_{t^*-1},v}$$ is maximized (minimized.)
\\[1pt]

\hrule
\noindent \rule{0pt}{5mm}

\begin{figure}
\begin{center}
\input{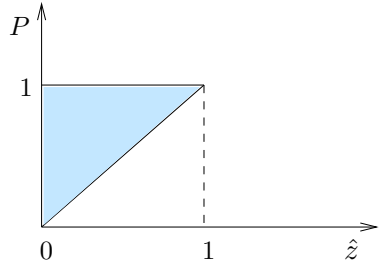}
\caption{Linear inequalities on $\hat z_{t;r_1,r_2,\ldots,r_{t-1},r_t}$ and $P_{t;r_1,r_2,\ldots,r_t}$.}
\label{why}
\end{center}
\end{figure}

\noindent The above system is not linear. By making the change of variables
\be
\lb{ChangeVariables}
\hat z_{t;r_1,r_2,\ldots,r_{t-1},r_t}=z_{t;r_1,r_2,\ldots,r_{t-1},r_t}P_{t;r_1,r_2,\ldots,r_t}
\ee
for $t=1,2,\ldots,t^*-1$,
we make it linear.
The domain for $\left(\hat z_{t;r_1,r_2,\ldots,r_t},P_{t;r_1,r_2,\ldots,r_t}\right)$ is represented on Figure \ref{why}.
Thus all the constraints as well as the objective function are linear.

\subsection {Computational Example}

\noindent The data on VerizonCM and Cisco call options (Table
\ref{data_options}) is taken from The Wall Street Journal of October
22 2002 to derive best possible bounds on the price of an option with a payoff function $(V+C-k)^+$, where $V$ is the price of the VerizonCM stock at the maturity, $C$ is the price of Cisco stock and $k$ is the strike price.  The prices of VerizonCM  and Cisco stocks on that day were \$37.75  and  \$11.22 correspondingly. 

Table \ref{exact approach} represents best possible bounds on the price of the option with maturity 52, if all given options are taken into account. Table \ref{exact approach1} gives bounds on the price of the same option, if only options of maturity 52 are taken into account. As expected, these bounds are looser. 

As can be noticed, bounds for strike prices smaller than 20 are extremely tight. That makes sense, since the probability that the sum of prices of the two given  stocks will go below 20, is almost zero. Thus, the option will be exercised almost surely and buying the option is almost equivalent to buying the two stocks.

\begin{table}
\centerline{\epsfig{file=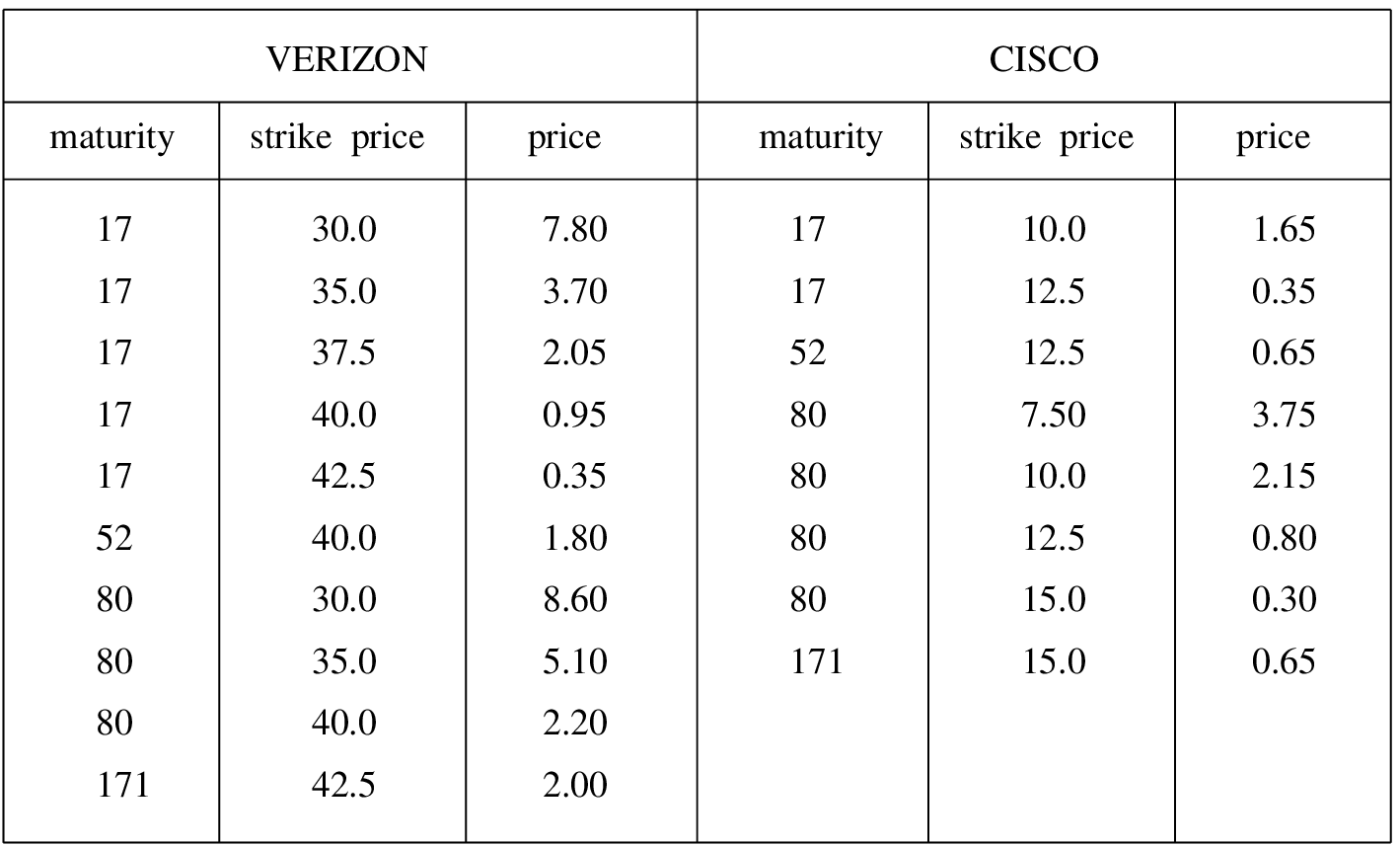, scale=0.82}}
\caption{\small{Data.}}
\label{data_options}
%

\centerline{\epsfig{file=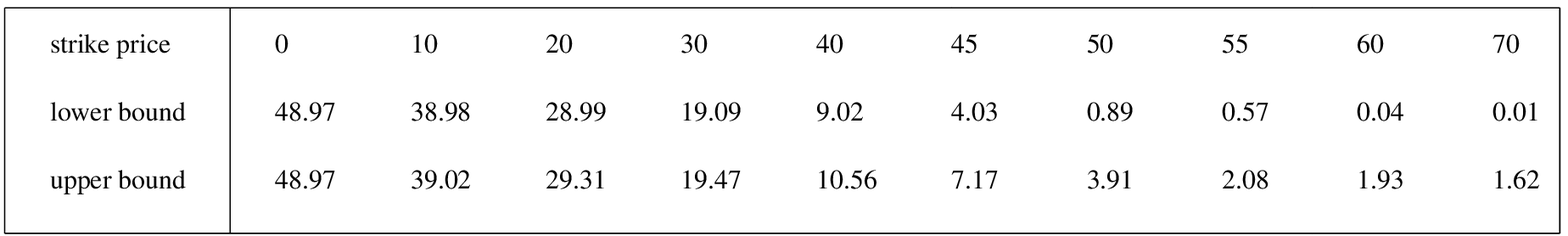, scale=0.82}}
\caption{\small{Best possible lower and upper bounds on the price of an option on the linear  combination of two stocks. Maturity is 52. Options of all maturities are taken into account.}}
\label{exact approach}

\centerline{\epsfig{file=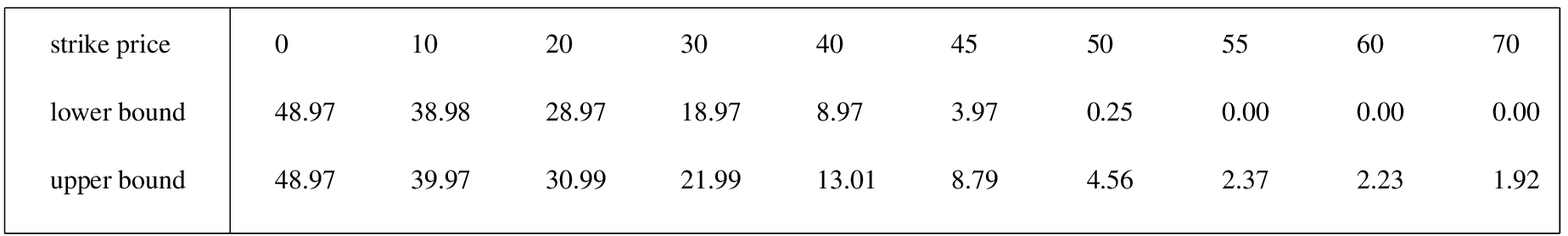, scale=0.82}}
\caption{\small{Best possible lower and upper bounds on the price of an option on the linear  combination of two stocks. Maturity is 52. Only options of this maturity are taken into account.}}
\label{exact approach1}

\end{table}

\subsection{Approximation approach}

In this section,  we develop an approximation approach which allows to
overcome the exponential growth (in $t^*$)
  of the number of variables.

Recall that we consider only martingales ${\bold X_t}$ such that $P({\bold X_t}\leq L)=1$ for some $L>0$. For now we will allow $g: \mathbb R^{2}_+\rightarrow\mathbb R$, to be any function continuous on $[0,L]\times[0,L]$ and develop a general approximation algorithm. Then we will show how the efficiency of the algorithm can be improved, if $g$ is a continuous piecewise linear function. Choose $\epsilon>0$ and designate by ${L}_{\epsilon}$ the set of all vertices of the $\epsilon$-square lattice in $[0,L]\times[0,L]\subset\mathbb R^2$.
\begin{theorem}
\label{one_dim}
Let ${\bold X_1},{\bold X_2},\ldots,{\bold X_{t^*}}$ be a two-dimensional Markov Martingale with values in $[0,L]\times[0,L]$. Then there exists a two-dimensional Markov Martingale ${\bold{\tilde X_1}},{\bold {\tilde X_2}},\ldots,{\bold{\tilde X_{t^*}}}$ such that for all $t=1,\ldots,t^*$, $P\left(\tilde {\bold X_{t}}\in L^t_{\epsilon}\right)=1$ and
$$P\left(\left|\tilde X_t^{{h}}-X_{t}^{{h}}\right|\leq t\epsilon\right)=1, \quad {h}= _{^I}, _{^{II}}.$$
Here $L^t_{\epsilon}$ is an $\epsilon$-square lattice in $[0,L+(t-1)\epsilon]\times  [0,L+(t-1)\epsilon]$.
\end{theorem}
\begin{proof}
Let ${\cal B}_1\subset{\cal B}_2\subset\cdots\subset{\cal B}_{t^*}$ be the family of $\sigma$-algebras corresponding to the martingale ${\bold X_1},{\bold X_2},\ldots,{\bold X_{t^*}}$.
Let ${\bold \Delta_t}:={\bold X_t}-{\bold X_{t-1}}$ for $t=2,3,\ldots,t^*$. Then
$${\bold X_t}={\bold X_1}+{\bold \Delta_2}+{\bold \Delta_3}+\cdots+{\bold \Delta_t}$$ and $E\left[{\bold \Delta_t}|{\cal B}_{t-1}\right]=0$.
\begin{figure}[h]
\centerline{\epsfig{file=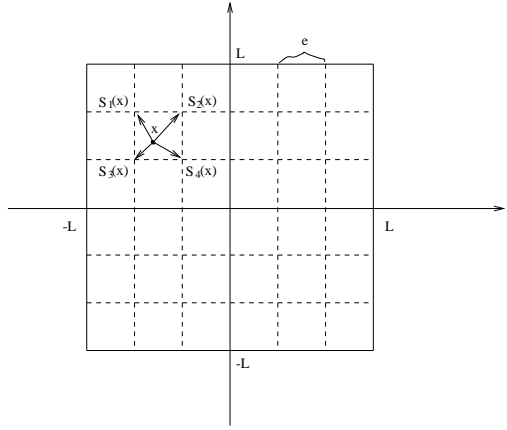}}
\caption{Lattice $\tilde L_{\epsilon}$.}
\label{lattice}
\end{figure}
\\[1pt]
Let $\tilde {L}_{\epsilon}$ designate the set of all vertices of the $\epsilon-$square lattice in $[-L,L]\times[-L,L]\in\mathbb R^2$ (Figure \ref{lattice}.) And let $S^j:\; [-L,L]\times[-L,L]\rightarrow \tilde L_{\epsilon}$, $j=1,2,3,4$,  be the family of functions, such that each function puts in correspondence to $x\in[-L,L]\times[-L,L]$ one of the closest to it four vertices of $\tilde L_{\epsilon}$.  

Notice that ${\bold \Delta_t}\in [-L,L]\times[-L,L]$ for all $t=2,\ldots,t^*$. Let us define a random variable ${\bold {\tilde\Delta_t}}$ such that conditioned on ${\bold \Delta_t}$, it can take only four values: $\left\{S^j({\bold \Delta_t)}-{\bold \Delta_t}\right.,$ $\left.j=1,2,3,4\right\}$ and such that 
$E\left[\left.{\bold {\tilde\Delta_t}}\right|{\cal B}_t\right]=0$. Moreover, we add a requirement, that conditioned on ${\cal B}_t$, ${\bold {\tilde{\Delta}_t}}$ is independent of any other random variables in consideration, including ${\bold{\tilde\Delta_1}},{\bold{\tilde\Delta_2}},\ldots,{\bold {\tilde\Delta_{t-1}}}$.
While ${\bold{\tilde\Delta_2}},{\bold{\tilde\Delta_3}},\ldots,{\bold{\tilde\Delta_{t^*}}}$ can be  defined based on ${\bold{\Delta_2}},{\bold \Delta_3},\ldots,{\bold \Delta_{t^*}}$ as described above, define ${\bold{\tilde\Delta_1}}$ in a similar way, but take ${\bold \Delta_1}={\bold X_1}$.

Then we have ${\bold \Delta_t}+{\bold {\tilde\Delta_t}}\in {{\tilde L_{\epsilon}}}$ and $P\left(\left|{\tilde\Delta_t}^{{h}}\right|\leq \epsilon\right)=1$, ${h}= _{^I},_{^{II}}$, $t=1,2\ldots,t^*$. 
Let $\sigma\left({\bold{\tilde\Delta_1}},{\bold{\tilde\Delta_2}},\ldots,{\bold{\tilde\Delta_{t}}}\right)$ denote the smallest $\sigma$-algebra for which ${\bold{\tilde\Delta_1}},{\bold{\tilde\Delta_2}},\ldots,{\bold{\tilde\Delta_{t}}}$ are all measurable and
$$
\tilde{\cal B}_t:=\sigma\left({\cal B}_t\cup \sigma\left({\bold{\tilde\Delta_1}},{\bold{\tilde\Delta_2}},\ldots,{\bold{\tilde\Delta_{t}}}\right)\right), t=1,2\ldots,t^*,$$
$$\begin{array}{lllllllll}
{\bold {\tilde X_1}}={\bold X_1}+{\bold{\tilde\Delta_1}},\\
{\bold{\tilde X_2}}={\bold X_2}+{\bold{\tilde\Delta_1}}+{\bold{\tilde\Delta_2}},\\
..........................................\\
{\bold {\tilde X_{t^*}}}={\bold X_{t^*}}+{\bold{\tilde\Delta_1}}+{\bold{\tilde\Delta_2}}+\cdots+{\bold{\tilde\Delta_{t^*}}}.
\end{array}
$$
\\[2pt]
${\bold {\tilde X_1}}\in L_{\epsilon}$ is obvious.
Now 
$${\bold {\tilde X_2}}={\bold X_2}+{\bold {\tilde\Delta_1}}+{\bold{\tilde\Delta_2}}={\bold X_1}+({\bold X_2}-{\bold X_1})+{\bold{\tilde\Delta_1}}+{\bold{\tilde\Delta_2}}=\left({\bold X_1}+{\bold {\tilde\Delta_1}}\right)+\left({\bold \Delta_2}+{\bold{\tilde\Delta_2}}\right)\in \tilde L_{\epsilon}.$$
But since also $P\left(\left|{{\tilde X_2}}^{{h}}-{X_2}^{{h}}\right|\leq 2\epsilon\right)=1$, ${h}= _{^I},_{^{II}}$, and ${\bold X_2}\in L_{\epsilon}\subset\tilde L_{\epsilon}$,  it must be that ${\bold {\tilde X_2}}\in L_{\epsilon}^2$.
Continuing recursively the theorem follows. 
\end{proof}
If ${\bold X_1}, {\bold X_2}, \ldots, {\bold X_{n}}$ is a two-dimensional martingale, that satisfies $E\left[f_{tj}^{{h}}({\bold X_t})\right]=C_{tj}^{{h}}$ for $(t,j,{h})\in{\cal T}$,  and
${\bold {\tilde X_1}}, {\bold{\tilde X_2}},\ldots, {\bold {\tilde X_n}}$ is a two-dimensional martingale, such that \\
$P\left(\left|{{\tilde X_t}}^{{h}}-{X_t}^{{h}}\right|\leq t\epsilon\right)=1$, ${h}= _{^I},_{^{II}}$, $t=1,\ldots,n$,
 then we must have 
\be
\lb{}
\left|E\left[f_{tj}^{{h}}\left({\bold {\tilde X_t}}\right)\right]-C_{tj}^{{h}}\right|\leq \epsilon t, \quad {h}= _{^I},_{^{II}},\quad t=1,\ldots,n 
\ee
and 
\be 
\lb{bound_g}
\left|E\left[g\left({\bold {\tilde X_{t^*}}}\right)\right]-E\left[g\left({\bold X_{t^*}}\right)\right]\right|\leq M\sqrt2\epsilon t^*
\ee
where $M$ is the smallest number,  such that 
$$\left|g({\bold x_2})-g({\bold x_1})\right|\leq M\left\Vert{\bold x_1}-{\bold x_2}\right\Vert_2$$
for all ${\bold x_1},{\bold x_2}\geq 0$,

\noindent Notice that if  $g(x,y)=|\alpha x+\beta y-k|^+$, then $\sqrt2M=\alpha+\beta$.
\\[3pt]
\begin{theorem} The following problem has its optimal solution converging to the optimal solution of problem\rf{conditions} as $\epsilon\rightarrow 0$:
\\[1pt]
Minimize (maximize) $E\left[g\left({\bold {\tilde X_{t^*}}}\right)\right]$ over all Markov martingales
${\bold {\tilde X_1}},{\bold {\tilde X_2}},\ldots,{\bold{\tilde X_{t^*}}}$, such that  $P\left({\bold {\tilde X_t}}\in  L_{\epsilon}^t\right)=1$ and
\label{converg}
\be
\lb{eps_1}
\left|E\left[f_{tj}^{{h}}\left({\bold {\tilde X_t}}\right)\right]-C_{tj}^{{h}}\right|\leq \epsilon t,\quad  \left\{({h},t,j)\in{\cal T}\;|\;t\leq t^*\right\},
\ee
\be
\lb{eps_2}
E\left[f_{fut,j}^{{h}}\left({\bold {\tilde X_{t^*}}}\right)\right]-C_{fut,j}^{{h}}\leq \epsilon t^*,\quad \left\{\left({h},_{^{fut}},j\right)\in{\cal T}\;|\;(k_{fut,j}^{{h}},C_{fut,j}^{{h}})\in {\cal F}^{{h}}\right\}.
\ee
\end{theorem}

\begin{proof}
Let us prove the theorem for the case when we want to maximize $E[g({\bold X_{t^*}})]$. 
For the minimization problem the proof is similar. 

Let $g^*$ be the optimal solution to the original problem (c) and $g^s$ the supremum as $\epsilon\rightarrow0$ of optimal solutions to problems with martingales' state spaces restricted to $L_{\epsilon}^t$ and subject to 
the constraints \rf{eps_1} and\rf{eps_2}.  
From Theorem \ref{one_dim} and \rf{bound_g} it immediately follows that $g^s\geq g^*$.
Suppose that there exists $\delta>0$, such that $g^s\geq g^*+\delta$. Then there exists a sequence of $\epsilon$'s, $\left\{\{\epsilon_n\}\;|\;\epsilon_1>\epsilon_2>\cdots\epsilon_n>\cdots;\; \lim_{n\rightarrow\infty}\epsilon_n\rightarrow0\right\}$, such that the limit of optimal solutions to problems corresponding to $\epsilon_1,\epsilon_2,\ldots,\epsilon_n,\ldots,$ is equal to $g^s$. Let $\left\{{\bold X^{\epsilon_n}_t}\right\}_{n\geq 1}$ be the sequence of martingales corresponding to the sequence $\{\epsilon_n\}$, such that  $\left\{{\bold X^{\epsilon_n}_t}\right\}_{n\geq 1}$  defines an optimal solution to the $\epsilon_n$-problem.

By definition of $\left\{{\bold X^{\epsilon_n}_t}\right\}$, it takes values in $L_{\epsilon_n}^{t^*}$ and
$$\left|E\left[f_{tj}\left({\bold X_t^{\epsilon_n}} \right)\right]-C_{tj}\right|\leq \epsilon_n t,\; \quad t=1,2,\ldots,t^*.$$
 Let $\mu^{\epsilon_n}$ designate the law of $\left\{{\bf X}^{\epsilon_n}_t\right\}_{t=1,\ldots,t^*}$. Then $\mu_{\epsilon_1},\mu_{\epsilon_2},\ldots,\mu_{\epsilon_n},\ldots$ have support in $[0,L+t^*\epsilon_1]^2$, which is a compact.
It follows that  the sequence of laws $\left\{\mu_{\epsilon_n}\right\}_{n\geq 1}$ is uniformly tight.
Consequently, there exists a subsequence $\mu_{\epsilon_{n({\mathfrak e})}}\rightarrow \mu$ for some law $\mu$ \cite[p.230]{C8}.  But then we have an admissible solution to our original problem with an optimal value greater than $g^*+\delta$, which is a contradiction. 
\end{proof}

\subsection{Algorithm for the approximation approach}

Let us enumerate all the nodes of $L_{\epsilon}^{t^*}$ and denote by ${\EuScript N}$ the set of all the indexes.
Designate by ${\bold Z}_{\mathfrak n}\in L_{\epsilon}^{t^*}$, the ${\mathfrak n}$-th node of $L_{\epsilon}^{t^*}$.
 Let $$P_{t;{\mathfrak n}}= P\left({\bold {\tilde X_t}}={\bold Z}_{\mathfrak n}\right)\quad \mbox{and}$$
$$P_{t;{\mathfrak n}_1,{\mathfrak n}_2}=P\left({\bold {\tilde X_t}}={\bold Z}_{{\mathfrak n}_1};\;{\bold {\tilde X_{t+1}}}={\bold Z}_{{\mathfrak n}_2}\right).$$
Notice that the distribution of $\tilde X_1,\tilde X_2,\ldots,\tilde X_{t^*}$ is uniquely determined by $P_{t;{\mathfrak n}_1,{\mathfrak n}_2}$, $t\in\{1,2,\ldots,t^*\}$, ${\mathfrak n}_1,{\mathfrak n}_2 \in{\EuScript N}$.

Following the same set of arguments as in the description of the algorithm for the exact solution to the one dimensional case, we can state that the problem consists in finding the collection of $P_{t;{\mathfrak n}}$ and $P_{t;{\mathfrak n}_1,{\mathfrak n}_2}$, $t\in\{1,2,\ldots,t^*\}$, ${\mathfrak n},{\mathfrak n}_1,{\mathfrak n}_2\in{\EuScript N}$, which solves the linear optimization problem:
$$
\begin{array}{rl}
\max~(\min) & \displaystyle \sum_{{{\mathfrak n}\in {\EuScript N}}}P_{t^*;{\mathfrak n}}\;g({\bold Z}_{\mathfrak n}) \vspace{3pt}\\
{\rm s.t.} &\sum_{{\mathfrak n}\in{\EuScript N}} P_{1;{\mathfrak n}}=1,\vspace{3pt}\\
& \displaystyle \sum_{{{\mathfrak n}\in {\EuScript N}}} P_{t;{\mathfrak n}_1,{\mathfrak n}}=P_{t;{\mathfrak n}_1}\quad\mbox{for all}\;\;t< t^*,\; {\mathfrak n}_1\in{\EuScript N},\vspace{3pt}\\
& \displaystyle \sum_{{{\mathfrak n}\in {\EuScript N}}} P_{t;{\mathfrak n},{\mathfrak n}_1}=P_{t+1;{\mathfrak n}_1}\quad\mbox{for all}\;\; t<t^*,\; {\mathfrak n}_1\in{\EuScript N},\vspace{3pt}\\
& \displaystyle \sum_{{{\mathfrak n}\in {\EuScript N}}}P_{t;{\mathfrak n}_1,{\mathfrak n}}{\bold Z}_{{\mathfrak n}}=P_{t;{\mathfrak n}_1}{\bold Z}_{{\mathfrak n}_1}\quad\mbox{for all}\; t<t^*,\; {\mathfrak n}_1\in{\EuScript N},\vspace{3pt}\\
& \displaystyle \;\sum_{{{\mathfrak n}\in {\EuScript N}}}P_{t;{\mathfrak n}}f_{tj}^{{h}}({\bold Z}_{\mathfrak n})-C_{tj}^{{h}}\leq t\epsilon,\quad \sum_{{{\mathfrak n}\in {\EuScript N}}}P_{t;{\mathfrak n}}f_{tj}^{{h}}({\bold Z}_{\mathfrak n})-C_{tj}^{{h}}\geq -t\epsilon\vspace{3pt}\\
& \displaystyle \mbox{for}\;\left\{({h},t,j)\in{\cal T}\;|\;t\leq t^*\right\},\vspace{3pt}\\
& \displaystyle \sum_{{{\mathfrak n}\in {\EuScript N}}}P_{t^*;{\mathfrak n}}f_{fut,j}^{{h}}({\bold Z}_{\mathfrak n})-C_{fut,j}^{{h}}\leq t^*\epsilon\\
& \displaystyle \mbox{for}\;\left\{({h},t,j)\in{\cal T}\;|\;t>t^*\right\},\\
& P_{t;{\mathfrak n}_1,{\mathfrak n}_2}\geq 0\quad\mbox{for all}\; t\leq t^*,\; {\mathfrak n}_1,{\mathfrak n}_2\in{\EuScript N}.
\end{array}
$$
As before this is a linear optimization problem.

\subsection{Approximation approach for the special case of a payoff function}
In this section, we  show how the efficiency of the approximation algorithm can be improved, if $g$ is a continuous piecewise linear function. 

Let $G^E = G_1^E\cup G_2^E\cup\cdots G_{t^*}^E$ and let graph $G$ be the graph with points of edges represented by $G^E$.
Then, if we restrict the state space of martingales to $G^E$, the solution to problem (c) will not change. More precisely:



\begin{proposition}
\label{th_graph}
Let ${\bold X_1},{\bold X_2},\ldots,{\bold X_{t^*}}$ be a nonnegative two dimensional Markov Martingale bounded by $L$. Then there exists a nonnegative Markov Martingale $\left\{{\bold {\tilde  X_1}},{\bold {\tilde X_2}},\ldots,{\bold {\tilde X_{t^*}}}\right\}$, with the state space in $G^E$, such that it  satisfies
$$E\left[f_{tj}^{{h}}({\bold X_t})\right]=E\left[f_{tj}^{{h}}\left({\bold {\tilde X_t}}\right)\right],\; t\leq t^*,\; ({h},t,j)\in {\cal T},$$ 
$$E\left[f_{fut,j}^{{h}}({\bold X_t})\right]=E\left[f_{fut,j}^{{h}}\left({\bold {\tilde X_t}}\right)\right],\; j=1,2,\ldots, U(_{fut},t^*),$$
$$E\left[g\left({\bold X_{t^*}}\right)\right]=E\left[g\left({\bold{\tilde X_{t^*}}}\right)\right].$$  
\end{proposition}

\begin{proof}
The proof is similar to  the proof of Theorem \ref{th_1} except that 
 we deal with two dimensional random variables and, thus, a two dimensional Brownian motion should be introduced instead of one dimensional.
\end{proof}

Let 
$$\Lambda^t_{\epsilon}:=\left\{(x,y)\in L^t_{\epsilon}\;|\;\exists\; (x_0,y_0)\in G_E: \left|x-x_0\right|\leq t\epsilon, \; |y-y_0|\leq t\epsilon\right\}.$$

Then Theorem \ref{converg} will hold if we take $\Lambda^t_{\epsilon}$ instead of $L^t_{\epsilon}$. The number of variables in this case grows only linear versus quadratically in the case of a general function $g$.

\subsection{Generalization to a multi-dimensional case}

The exact solution, the approximation approach for an arbitrary payoff function $g$ and the approximation approach for a continuous piecewise linear payoff function $g$, can all be extended from a two dimensional case to a multiple dimensional case. All theorems and definitions of this chapter can be reformulated for an $n$-dimensional case, and proofs will be identical up to a dimensionality.
However the number of variables in the linear optimization 
problem will grow exponentially in the  dimension.

\subsection{Generalization to the case of options with continuous piecewise linear payoff functions}
Across the paper we considered the case when prices of simple call options are given and bounds on prices of options with continuous piecewise linear functions must be found.  From financial point of view this formulation of the problem can be motivated by the fact that simple call and put options are more liquid and prices on them are readily available in the market, compared to more complex exotic options. Notice, that conditions on European put options can be easily expressed as conditions on call options.  

However, the problem could be extended to a more general case when options with continuous piecewise linear functions are given and bounds on similar type of options are to be determined. If there were no options of maturities larger than $t^*$, then all that has to be changed in the solution is definitions of the set ${\cal K}$ in the one dimensional case and graphs $G_t$ in the two dimensional case. In particular, in the one dimensional case,  in section 3.2, ${\cal K}$ must include all points of $\mathbb R^+$  where payoff functions of given and target options change their derivatives. Similarly, for the two dimensional case, in section 4.2, lines defined in (a) and (b) of Definition {\ref{Definition_Graph_G}} by strike prices of call options would have to be defined by points of $\mathbb R^+$ where payoff functions of options on each asset change their derivatives. Then theorems that refer to definitions of set ${\cal K}$ and graphs $G_t$  would remain valid.

However, the simplification that was made for treatment of future conditions (sections 3.1 and 4.1) is not possible if payoff functions of options of maturities larger than $t^*$ are piecewise linear, since these functions are not necessarily convex. As a result, the martingale has to be constructed not only up to time $t^*$, but up to the time of the maximum options maturity, and future and present constraints have to be handled in the same way as past constraints. In the one dimensional case, since the state space of the martingale at each time is the same, there is no principal difference. All that one  has to do  in addition to changes in the definition of set ${\cal K}$ is to extend the set of unknowns  $P_{t;{\mathfrak n}}$ and $P_{t;{\mathfrak n}_1,{\mathfrak n}_2}$ up to the maximum time $t$ and change inequalities to equalities in the future constraints. The problem remains a linear optimization problem.
In the two dimensional case, however, there would not necessarily exist a martingale that takes values at the vertices of graph $G_{t^*}$ at time $t^*$ and satisfy all future constraints. As was shown in the proof of Theorem \ref{th_graph2}, there does exist  a martingale that satisfies all given constraints and at time $t^*$ has an atomic conditional distribution with at most one atom per each region of $G_{t^*}^R$, conditioned on the sequence of all prior visited regions.  However, since regions of $G_{t^*}^R$ are not necessarily rectangular, the simple form of parametrization\rf{parametrization} and the change of variables\rf{ChangeVariables} can not be applied. Therefore, we don't think it is possible to reduce present constraints and the objective function to the linear form in this case.

It is easy to see that the approximation approach for the two dimensional case can be applied to the solution of the problem, when prices of options with continuous piecewise linear payoff functions are given. In this case, one only needs to extend the set of unknowns  $P_{t;{\mathfrak n}}$ and $P_{t;{\mathfrak n}_1,{\mathfrak n}_2}$ from time $t^*$ to the maximum maturity and  use  equalities instead of inequalities for future constraints.


\section{Static arbitrage bounds on basket option prices}
In this section, we address Problem (c) defined as follows:
\\[1pt]

\noindent \parbox[c]{13.5cm}{\noindent {\bf  Problem (d)}
\hrule
\noindent \rule{0pt}{5mm} 
\noindent Given ${\bold C}\in \mathbb R^{m}_+,\;{\bold k}\in\mathbb R^{m}_+,\; \omega_i\in\mathbb R^n,\; i=1,\ldots,m$ and $k_0>0$, $w_0\in\mathbb R^{m}_+$, find the upper and lower bounds on 
\be
\lb{funct_1}
E\left({\bold \omega}_0^{T}{\bold X}-{k}_0\right)^+,
\ee
with respect to distributions of an $n-$dimensional random variable ${\bold X}$ with finite expectation and support in $\mathbb R^{n}_+$ under the condition 
\be
\lb{funct_2}
E\left(\omega_i^T{\bold X}-{k}_i\right)^+={C}_i, \quad i=1,\ldots,m.
\ee
\hrule
\noindent \rule{0pt}{5mm}
}

\noindent To solve the problem we  make an additional assumption that the support of ${\bold X}$ is bounded, that is, there exists $L>0$, such that $P\left(X^{{h}}\leq L\right)=1$ for any ${h}= 1,2,\ldots,n$.

Let $f_i:\mathbb R^{n}_+\rightarrow {\mathbb R}_+$ be the family of measurable functions defined as
$$f_i({\bold X})=\left({\bold w}_i^T{\bold X}-{k}_i\right)^+, \quad i=0,1,\ldots,m.$$
Let $\cup_j D_j=[0,L]^n$ be the partition of $[0,L]^n$, into  subsets $D_j$ of $\mathbb R^{n}_+$ such that in the interior of each subset $D_j$ all functions $f_i$, $i=0,1,\ldots,m$ are affine.

%

Let $G\subset[0,L]^n$ 
be the graph in $\mathbb R^{n}_+$ 
formed by intersecting hyperplanes  
\be
\lb{hyperplanes}
\omega^T_i{\bold X}=k_i, \quad i=1,2,\ldots,m,\nonumber\\
X^{{h}}=0, \quad {h}=1,2,\ldots,n, \\
X^{{h}}=L, \quad {h}=1,2,\ldots,n.\nonumber
\ee
Let us enumerate these hyperplanes and call them $h_1, h_2,\ldots, h_{2n+m}$ correspondingly. 
We call a vector ${\bold r}\in \mathbb [0,L]^n$ a vertex of $G$ if there are $n$ independent hyperplanes of\rf{hyperplanes} which  intersect at ${\bold r}$. Let $G_V$ designate the set of all vertices. Notice, that this set contains no more than $C_{2n+m}^n=\frac{(2n+m)!}{(n+m)! n!}$ elements. 
\begin{corollary}
For each random variable ${\bold X}$ with support in $[0,L]^n$ there exists a random variable ${\bold {\bar X}}$ such that $P\left({\bold {\bar X}}\in G_V\right)=1$ and $E\left[f_i({\bold X})\right]=E\left[f_i({\bold {\bar X}})\right]$, $i=0,1,\ldots,m$.
\end{corollary}

Thus, to find optimal bounds we  solve a linear optimization
 problem with unknown variables being weights of the atomic distribution. The atoms are located in the vertices of the defined graph.

\end{document}